\documentclass[10pt]{article}
\usepackage{amsmath,amssymb,amsthm,amscd}
\numberwithin{equation}{section}

\DeclareMathOperator{\rank}{rank}
\DeclareMathOperator{\sing}{Sing}
\DeclareMathOperator{\Aut}{Aut}
\DeclareMathOperator{\spec}{Spec}

\DeclareMathOperator{\Ext}{Ext}

\DeclareMathOperator{\aut}{Aut}

\DeclareMathOperator{\image}{Im}
\DeclareMathOperator{\Int}{Int}

\DeclareMathOperator{\ob}{Ob}

\DeclareMathOperator{\supp}{Supp}
\DeclareMathOperator{\coker}{Coker}

\def\ext{\mathop{\cE xt\skp}}

\def\A{{\mathfrak A}}

\def\C{{\mathfrak C}}
\def\W{{\mathfrak W}}

\def\g{\mathfrak g}

\def\h{{\mathfrak h}}
\def\M{{\mathfrak M}}

\def\R{{\mathfrak R}}
\def\fS{{\mathfrak S}}
\def\W{{\mathfrak W}}

\def\cA{{\mathcal A}}
\def\cC{{\mathcal C}}
\def\cB{{\mathcal B}}

\def\cD{\mathcal D}
\def\cE{{\mathcal E}}
\def\cF{{\mathcal F}}
\def\cG{{\mathcal G}}
\def\cH{{\mathcal H}}
\def\cI{{\mathcal I}}

\def\cK{{\mathcal K}}
\def\cL{{\mathcal L}}
\def\cM{{\mathcal M}}
\def\cN{{\mathcal N}}
\def\cO{{\mathcal O}}
\def\cP{{\mathcal P}}
\def\cQ{{\mathcal Q}}

\def\cS{{\mathcal S}}
\def\cT{{\mathcal T}}
\def\cU{{\mathcal U}}
\def\cV{{\mathcal V}}
\def\cW{{\mathcal W}}
\def\cX{{\mathcal X}}
\def\cY{{\mathcal Y}}

\def\bB{{\mathbf B}}
\def\bC{{\mathbf C}}
\def\bE{{\mathbf E}}

\def\bv{{\mathbf v}}

\def\bh{{\mathbf h}}

\def\bg{{\mathbf g}}
\def\bphi{{\boldsymbol \phi}}

\def\ZZ{{\mathbb Z}}
\def\QQ{{\mathbb Q}}
\def\RR{{\mathbb R}}
\def\CC{{\mathbb C}}

\newtheorem{prop}{Proposition}[section]
\newtheorem{theo}[prop]{Theorem}
\newtheorem{lemm}[prop]{Lemma}
\newtheorem{coro}[prop]{Corollary}

\newtheorem{defi}[prop]{Definition}

\def\begeq{\begin{equation}}
\def\endeq{\end{equation}}

\def\and{\quad{\rm and}\quad}

\def\bl{\bigl(}
\def\br{\bigr)}
\def\lbe{_{\beta}}

\def\bs{{\mathbf s}}

\def\clo{\mathop{\rm cl\skp}}

\def\dual{^{\vee}}

\def\dbar{\bar\partial}

\def\ddotw{{\mathcal D}_{\tilde w}^{\bullet}}
\def\dbul{{\mathcal D}^{\bullet}}

\def\etilu{\bE_{\tilde U}}
\def\eqdef{: =}
\def\etil{\tilde{\bB}}
\def\eps{\epsilon}
\def\ueps{^{\epsilon}}
\def\extfc{\Ext^1(\ddotw,\cO\lsig)\shar}
\def\ebul{\cE^{\bullet}}

\def\gtilu{G_{\tilde U}}

\def\hf{{1/2}}
\def\uhf{^{\hf}}

\def\lsa{_{S,\alpha}}
\def\lalpbe{_{\alpha\beta}}
\def\lsat{_{S,\alpha,2}}
\def\lst{_{S,2}}
\let\lra=\longrightarrow
\def\lsta{_{\ast}}
\def\lsig{_{\Sigma}}
\def\lbe{_{\beta}}
\def\lex{_{{\rm ex}}}
\def\lalp{_{\alpha}}
\def\lamzo{\Lambda^{0,1}}
\def\ltilu{_{\tilu}}
\def\ltilv{_{\tilde V}}
\def\lrv{_{R_V}}
\def\ltilw{_{\tilde W}}
\def\lkl{_{k,l}}

\def\lexp{_{\rm exp }}
\def\lcpt{_{\rm cpt}}

\def\mgna{\M_{g,n}(X,A)}
\def\mapright#1{\,\smash{\mathop{\lra}\limits^{#1}}\,}
\def\mapto{\xrightarrow}
\def\mh{\!:\!}
\def\mgn{\M_{g,n}}
\def\mgnk{\M_{g,n+k}}
\def\mgnab{\overline{\M}_{g,n}(X,A)}
\def\mgnb{\overline{\M}_{g,n}}
\def\om{\Omega}
\def\mm{{\mathfrak m}}

\def\Om{\Omega}

\def\Phitil{\tilde{\Phi}}
\def\pri{^{\prime}}
\def\pitilu{\pi_U}
\def\phitilu{\Phi\ltilu}

\def\sub{\subset}

\def\sta{^{\ast}}
\def\sha{^{!}}
\def\subb{\Subset}
\def\skp{\hspace{1pt}}
\def\shar{^{\dag}}

\def\tilu{\tilde U}
\def\tilw{\tilde W}
\def\tilv{\tilde V}
\def\tile{\tilde{{\mathbf E}}}

\def\top{^{\rm top}}

\def\upmo{^{-1}}
\def\ua{^{\cA}}
\def\uj{^{J}}
\def\util{\tilde U}
\def\ulk{^l_k}
\def\uan{^{\text{an}}}
\def\ualg{^{\text{alg}}}

\def\upre{^{{\rm pre}}}

\def\vir{^{{\rm vir}}}

\def\lab{\ }
\def\lab{\label}

\title{Comparison of the algebraic and
the symplectic
Gromov-Witten invariants}

\author{Jun Li\thanks{Supported partially by NSF grant, Sloan fellowship
and Terman fellowship}
\\Department of Mathematics\\
Stanford University\\
Stanford, CA 94305 \\ and \\
Gang Tian\thanks{Supported partially by 
NSF grants}\\Department of Mathematics
\\Massachusetts Institute of Technology\\
Cambridge, MA 02139}
\date{}

\begin{document}

\maketitle

\section{Introduction}

As Witten suggested in [W1], [W2], the GW-invariants for
a symplectic manifold $X$ are multi-linear maps
\begin{equation}
\gamma_{A,g,n}^X: H\sta(X;\QQ)^{\times n}\times H\sta(\mgnb;\QQ)
\lra \QQ,
\lab{eq:0.1}
\end{equation}
where $A\in H_2(X,\ZZ)$ is any homology class,
$n$, $g$ are two non-negative integers, and $\mgnb$ is the Deligne-Mumford compactification of $\mgn$, the space of smooth $n$-pointed genus $g$ curves.
The basic idea of defining these invariants is to enumerate 
holomorphic maps from
Riemann surfaces to the manifolds. To illustrate this, we let
$X$ be a smooth projective manifold and form the moduli space $\mgna$
of all holomorphic maps $f\mh\Sigma\to X$ from smooth $n$-pointed
Riemann surfaces $(\Sigma;x_1,\ldots,x_n)$ to $X$ such that 
$f\lsta([\Sigma])=A$. $\mgna$ is a quasi-projective scheme and its
expected dimension can be calculated using the Riemann-Roch theorem. We 
will further elaborate the notion of expected dimension later, and for the
moment we will denote it by $r\lexp$. 
Note that it depends implicitly on the choice of
$X$, $A$, $g$ and $n$. When $r\lexp=0$, then $\mgna$ is expected to
be discrete. If $\mgna$ is discrete, then the degree of $\mgna$, considered as
a $0$-cycle, is a GW-invariant of $X$. We remark that we have and will 
ignore the issue of non-trivial automorphism groups of maps in $\mgna$
in the introduction.
When $r\lexp>0$, then $\mgna$ is expected to have pure dimension $r\lexp$. If
it does, then we pick $n$ subvarieties of $X$, say $V_1,\ldots,V_n$, so that
their total codimension is $r\lexp$. We then form a subscheme of 
$\mgna$ consisting of maps $f$ so that $f(x_i)\in V_i$. This subscheme
is expected to be discrete. It it does, then its degree is the GW-invariant
of $X$. Put them together, we can define the GW-invariants 
$\gamma_{A,g,n}^X$ of $X$. This is similar to
construction of the Donaldson polynomial invariants for 4-manifolds.

Here are the two big {\sl ifs} in carry out this program are

{\bf Question I}: Whether the moduli scheme $\mgna$ has pure dimension $r\lexp$.

{\bf Question II}: Whether the subschemes of $\mgna$ that satisfy certain
incidence relations have the expected dimensions.

Similar to Donaldson polynomial invariants, the affirmative answer 
to the above two questions are in general not guaranteed. 
One approach to overcome this difficulty, 
beginning with Donaldson's invariants of 4-manifolds,
is to ``deform'' the moduli problems and hope that the answers to the 
``deformed'' moduli problems are affirmative. In the case of GW-invariants, one
can deform the complex structure of the smooth variety $X$ to not
necessary integrable almost complex structure $J$ and study the same
moduli problem by replacing holomorphic maps with pseudo-holomorphic
maps. This was investigated by Gromov in [Gr], Ruan [Ru], 
in which he constructed certain GW-invariants of rational type for 
semi-positive symplectic manifolds. The first mathematical theory 
of GW-invariants came from the work of Ruan and the second author, 
in which they found that the right set up of GW-invariants for 
semi-positive manifolds can be provided by using the moduli of 
maps satisfying non-homogeneous Cauchy-Riemann equations.
In this set up, they constructed the GW-invariants of all semi-positive
symplectic manifolds and proved fundemantal properties of
these invariants. All Fano-manifolds and Calabi-Yau
manifolds are special examples of semi-positive symplectic
manifolds. Also any symplectic manifold of complex dimension less than $4$
is semi-positive.

Attempts to push this to cover general symplectic manifolds so
far have failed. New approaches are needed in order
to get a hold on the
GW-invariants of general varieties (or symplectic manifolds).  
The first step
is to  convert the problem of counting mappings, which essentially is homology
in nature, into the frame work of cohomology theory of the moduli problem.
More precisely, we first compactify the moduli space $\mgna$ to, say,
$\mgnab$. We require that the obvious evaluation map 
$$e: \mgna\lra X^n
$$
that sends $(f;\Sigma;x_1,\ldots,x_n)$ to $(f(x_1),\ldots,f(x_n))$ extends to 
$$\bar e: \mgnab\lra X^n.
$$
We further require that if $\mgna$ has pure dimension $r\lexp$, then
$\mgnab$ supports a fundamental class
$$[\mgnab]\in
H_{2r\lexp}(\mgnab;\QQ).
$$
Then the GW-invariants of $X$ are multi-linear
maps
\begin{equation}
\gamma^X_{A,g,n}: H\sta(X)^{\times n}\times H\sta(\mgnb)\lra \QQ
\lab{eq:0.4}
\end{equation}
that send $(\alpha,\beta)$ to
$$\gamma^X_{A,g,n}(\alpha,\beta)=\int_{[\mgnab]}{\bar e}\sta(\alpha)
\cup\pi\sta(\beta).
$$
where $\pi\mh\mgnab\to\mgnb$ is the forgetful map. Note that in such cases 
the GW-invariants are defined without reference to the answer to
{\sl question II}. 

Even when the answer to {\sl question I} is negative, we can still define
the GW-invariants if a virtual moduli cycle
$$[\mgnab]\vir\in H_{2r\lexp}(\mgnab;\QQ)
$$
can be found that function as the fundamental cycle
$[\mgnab]$ should the dimension of $\mgna$
is $r\lexp$. In this case, we simply
define 
$\gamma^X_{A,g,n}$ as before with $[\mgnab]$ replaced by $[\mgnab]\vir$.

The standard compactification of $\mgna$ is the moduli space of stable
morphisms from $n$-pointed genus $g$ curves, possibly nodal, to $X$
of the prescribed fundamental class. This was first studied for 
pseudo-holomorphic maps by T. Parker and J. Wolfson \cite{PW} and 
in algebraic geometry
by Kontsevich \cite{Ko}.
Because points of the compactification $\mgnab$ are maps $f$ whose domains
have $n$-marked points $x_1,\ldots,x_n$, the evaluation map 
$e$ extends canonically to $\bar
e$ that sends such map $f$ to $(f(x_1),\ldots,f(x_n))$.

The virtual moduli cycles $[\mgnab]\vir$ for projective variety $X$ were first
constructed by the authors. Their idea is to construct a virtual
normal cone embedded in a vector bundle based on the obstruction theory of
stable morphisms \cite{LT1}. An alternative construction of such cones
was achieved by Behrend and Fantechi \cite{BF,
Be}. For general symplectic manifolds, such virtual moduli cycles were 
constructed by the authors,
and independently, by Fukaya and Ono \cite{FO, LT2}. Shortly after
them, B. Siebert [Si] and later, Y. Ruan [Ru2]
gave different constructions of such virtual moduli cycles.
Both Siebert and Ruan's approach needs to construct
global, finite-dimensional resolutions of so called cokernel bundles
(cf. \cite{Si} and \cite{Ru2}, Appendix).

However, one question remains to be investigated. Namely, if $X$ is
a smooth projective variety then on one hand we have the algebraically
constructed GW-invariants, and on the other hand, by viewing $X$ as a 
symplectic manifold using the K\"ahler form on $X$, we have the
GW-invariants constructed using analytic method. These two approaches
are drastically different. One may expect, although far from clear, that 
for smooth projective varieties the 
algebraic GW-invariants and their symplectic counterparts are identical.

The main goal of this paper is to prove what was expected is indeed true.

\begin{theo}
Let $X$ be any smooth projective variety with a K\"ahler form $\omega$.
Then the algebraically constructed GW-invariants of $X$ coincide with
the analytically constructed GW-invariants of the symplectic
manifold $(X^{{\rm top}},\omega)$.
\end{theo}

This result was first announced in [LT2].
Its proof was outlined in [LT3]. 
During the preparation of the paper, we learned 
from B. Siebert that he was able to prove a similar result.

We now outline the proof of our Comparison Theorem. We begin with a few
words on the algebraic construction of the virtual moduli cycle.
Let $w\in\mgnab$ be any point associated to the stable morphism 
$f\mh \Sigma\to X$. It follows
from the deformation theory of stable morphisms that there is a complex
$\cC_w$, canonical up to quasi-isomorphisms, such that its first cohomology
$\cH^1(\cC_w)$ is the space of the first order deformations of the map $w$,
and its second cohomology $\cH^2(\cC_w)$ is the
obstruction space to deformations of the map $w$. 
Let $\varphi_w$ be a Kuranishi map of the obstruction theory of $w$.
Note that $\varphi_w$ is the germ of a holomorphic map from
a neighborhood of the origin $o\in\CC^{m_1}$ to $\CC^{m_2}$,
where $m_i=\dim\cH^i(\cC_w)$.
Let $\hat o$ be the formal completion of $\CC^{m_1}$ along $o$
and let $\hat w$ be the subscheme of $\hat o$ defined by the 
vanishing of $\varphi_w$. Note that $\hat w$ is isomorphic to the 
formal completion of $\mgnab$ along $w$
(Here as before we will ignore the issue of non-trivial automorphism 
groups of maps in $\mgnab$).
This says that ``near'' $w$, the scheme $\mgnab$
is a ``subset'' of $\CC^{m_1}$ defined by the vanishing of $m_2$-equations.
Henceforth, it these equations are in general position, them
$\dim\hat w=m_1-m_2$,
which is the expected dimension $r\lexp$ we mentioned before. The case where
$\mgnab$ has dimension bigger than $r\lexp$ is exactly when the
vanishing locus of these $m_2$- equations  in
$\varphi_w$ do not meet properly near $o$.
Following the excess intersection theory of Fulton and MacPherson \cite{Fu},
the ``correct'' cycle should come from first constructing the normal cone 
$C_{\hat w/\hat o}$ to $\hat w$ in $\hat o$, which is canonically a subcone
of $\hat w\times \CC^{m_2}$, and then intersect the cone with
the zero section of $\hat w\times\CC^{m_2}\to\hat w$. 
The next step is to patch these cones
together to form a global cone over $\mgna$. The main
difficulty in doing so comes from the fact that the
dimensions
$\cH^2(\cC_w)$ can and do vary as $w$ vary, only
$\dim\cH^1(\cC_w)-\dim\cH^2(\cC_w)$ is a topological number.
This makes the cones
$C_{\hat w/\hat o}$ to sit inside bundles of varying ranks. To overcome this
difficulty, the authors came with the idea of finding a global $\QQ$-vector
bundle $E_2$ over $\mgnab$ and a subcone
$N$ of $E_2$ such that near fibers over $w$, 
the cone $N$ is a fattening of the cone
$C_{\hat w/\hat o}$ (See section 3 or \cite{LT1} for more details).
In the end, we let $j$ be the zero section of $E_2$ and let $j\sta$
be the Gysin map 
$$ A\lsta E_2\lra A\lsta\mgnab,
$$
where $A\lsta$ denote the Chow-cohomology group (see \cite{Fu}). 
Then the algebraic virtual moduli cycle is
$$[\mgnab]\vir=j\sta([N])\in A_{r\lexp}\mgnab.
$$

Now let us recall briefly
the analytic construction of GW-invariants of
symplectic manifolds. Let $(X,\omega)$ be any smooth symplectic manifold
with $J$ a tamed almost complex structure. For $A$, $g$ and $n$ as before, we can form
the moduli space of $J$-holomorphic maps $f\mh \Sigma\to X$ where $\Sigma$
are $n$-pointed smooth Riemann surfaces such that $f\lsta([\Sigma])=A$. 
We denote this space by $\mgna\uj$. 
It is a finite dimensional topological space. 
As before, we compactify it to include all
$J$-holomorphic maps whose domains are possibly with nodal singularities.
We denote the compactified space by $\mgnab\uj$. To proceed, we will embed $\mgnab\uj$
inside an ambient space $\bB$ and realize it as the vanishing locus of a section of a
``vector bundle''.
Without being precise, the space $\bB$ is the space of all {\sl smooth} maps $f
\in\bB$ from
possibly nodal $n$-pointed Riemann surfaces to $X$, the fiber of the bundle over $f$
are all $(0,1)$-forms over  domain$(f)$ with values in $f\sta T_X$ and the section is the
one that sends $f$ to $\dbar f$. 
We denote this bundle by $\bE$ and the section by $\Phi$.  
Clearly, $\Phi\upmo(0)$ is homeomorphic to $\mgnab\uj$.  
Defining the GW-invariants of $(X,\omega)$ is essentially about constructing the Euler
class of $[\Phi\mh\bB\to\bE]$. 
This does not make much sense since $\bB$ is an
infinite dimensional topological space.
Although at each $w\in\Phi\upmo(0)$ the formal differential $d\Phi(w)\mh T_w\bB\to\bE_w$
is Fredholm, which has real index $2r\lexp$, the conventional perturbation scheme does not
apply directly since near maps in $\bB$ 
whose domains are singular the space $\bB$ is not smooth and
$\bE$ does not admit local  trivializations. 
To overcome this difficulty, the authors introduced the notion of 
weakly $\QQ$-Fredholm
bundles, and showed in \cite{LT2} that $[\Phi\mh \bB\to\bE]$ is a weakly $\QQ$-Fredholm
bundle and that any weakly $\QQ$-Fredholm bundle admits an Euler class. 
Let 
$$
e[\Phi\mh\bB\to\bE]\in H_{2r\lexp}(\bB;\QQ)
$$
be the the Euler class of $[\Phi\mh \bB\to\bE]$. 
Since the evaluation map of $\mgna\uj$
extends to an evaluation map ${\mathbf e}\mh\bB\to X^n$, 
the Euler class, which will
also be referred to as the symplectic virtual cycle of $\mgnab\uj$, defines a
multi-linear map
$\gamma_{A,g,n}^{X,J}$ as in \eqref{eq:0.1}. 
We will review the notion of weakly smooth Fredholm bundles in section 2. Here
to say the least, $[\Phi\mh\bB\to\bE]$ is weakly Fredholm means that near each point of
$\Phi\upmo(0)$ we can find a finite rank subbundle $V$ of $\bE$ such that
$W=\Phi\upmo(V)$ is a smooth finite dimensional manifold, $V|_W$ is a smooth vector
bundle and the lift $\phi\mh W\to V|_W$ of $\Phi$ is smooth. 
(Note that the rank of $V$ may vary but $\dim_{\RR}W-\rank_{\RR} V=2r\lexp$). 
For such finite models
$[\phi\mh W\to V|_W]$, which are called weakly smooth approximations, we can perturb
$\phi$ slightly to obtain
$\phi\pri$ so that $\phi^{\prime-1}(0)$ are smooth manifolds in $W$. To construct the
Euler class, we first cover a neighborhood of $\Phi\upmo(0)$ in $\bB$ by finitely many
such approximations that satisfy certain compatibility condition. We then perturb
each section in the approximation and obtain a collection of locally closed
$\QQ$-submanifolds of $\bB$ of dimension $2r\lexp$. By imposing certain compatibility
condition on the perturbations, this collection of $\QQ$-submanfolds patch
together to form a $2r\lexp$-dimensional cycle in $\bB$, which represents a homology
class in
$H_{2r\lexp}(\bB;\QQ)$. This is the Euler class of $[\Phi\mh\bB\to\bE]$.

Now we assume that 
$X$ is a smooth projective variety and $\omega$ is a K\"ahler 
form of $X$. Let $J$ be the complex structure of $X$. Then $\mgnab$ is
homeomorphic to $\mgnab\uj$. Hence the two GW-invariants $\gamma_{A,g,n}^X$
and $\gamma_{A,g,n}^{X,J}$ are identical if the homology classes
$[\mgnab]\vir$ and $e[\Phi\mh\bB\to\bE]$ will be
identical. Here we view $[\mgnab]\vir$ as a class in 
$H\lsta(\bB;\QQ)$ using 
$$\mgnab\sim_{{\rm homeo}}\mgnab\uj\sub\bB.
$$
To illustrate why these two classes are equal, let us first look at 
the following simple
model. Let $Z$ be a compact smooth variety and let $E$ be a holomorphic
vector bundle over $Z$ with a holomorphic section $s$. 
There are two ways to construct the Euler classes of
$E$. One is to perturb $s$ to a smooth section $r$ so that the graph of $r$ is
transversal to the zero section of $E$, and then define the Euler class of 
$E$ to be the
homology class in $H\lsta(Z;\QQ)$ of $r\upmo(0)$. This is
the topological construction of the Euler class of $E$. The algebraic
construction is as follows. Let $t$ be a large scalar and let 
$\Gamma_{ts}$ be the graph of
$ts$ in the total space of $E$. Since $s$ is an algebraic section, it follows
that the limit 
$$\Gamma_{\infty s} =\lim_{t\to\infty}\Gamma_{ts}
$$
is a complex dimension $\dim Z$ cycle supported on union of
subvarieties of
$E$. We then let $r$ be a smooth section of $E$
in general position and let $\Gamma_{\infty s}\cap
\Gamma_r$ be their intersection. Its image in $Z$ defines a
homology class, which is the image of the Gysin map $j\sta([C])$, where $j$
is the zero section of $E$. The reason that 
$$e(E)=[r\upmo(0)]=j\sta([\Gamma_{\infty s}])\in H\lsta(Z;\QQ)
$$
is that if we choose $r$ to be in 
general position, then 
$$[r\upmo(0)]=[\Gamma_r\cap\Gamma_0]=[\Gamma_r\cap\Gamma_s],
$$
and the family $\{\Gamma_{ts}\cap\Gamma_r\}_{t\in[1,\infty]}$ 
forms a homotopy of the cycles $\Gamma_s\cap\Gamma_r$ and 
$\Gamma_{\infty s}\cap \Gamma_r$. One important remark is that
the cone $\Gamma_{\infty s}$ is contained in $E|_{s\upmo(0)}$ and 
the intersection of $\Gamma_{\infty s}$ and $\Gamma_r$ in
$E$ is the same as their
intersection in $E|_{s\upmo(0)}$.

Back to our construction of GW-invariants, the analytic construction of
GW-invariants, which was based on perturbations of sections 
in the finite models
(weakly smooth approximations) $[\phi\mh W\to V|_W]$, is clearly a
generalization of the topological construction of the Euler classes
of vector bundles.
As to the algebraic construction of GW-invariants, it is based on a
cone in a $\QQ$-vector bundle over $\mgnab$. Comparing to the algebraic
construction of the Euler class of $E\to Z$, what is missing is the
section $s$  and that the cone is the limit of the graphs of the
dilations of $s$. Following \cite{LT1}, the cone $\Gamma_{\infty s}$
only relies on the  restriction of $s$ to an ``infinitesimal'' neighborhood of
$s\upmo(0)$ in $Z$, and can also be reconstructed using the Kuranishi maps
of the obstruction theory to deformations of points in $s\upmo(0)$ induced
by the defining equation $s=0$. Along this line, to each finite model
$[\phi\mh W\to V|_W]$ we can form a cone $\Gamma_{\infty \phi}=\lim
\Gamma_{t\phi}$ in $V|_{\phi\upmo(0)}$. Hence to show that the two virtual
moduli cycles coincide, it suffices to establish a relation, similar to
quasi-isomorphism of complexes, between the cone
$N$ constructed based the obstruction theory of $\mgnab$ and the collection
$\{\Gamma_{\infty\phi}\}$. In the end, this
is reduced to showing that the
obstruction theory to deformations of maps in $\mgnab$ is identical to the
obstruction theory to deformations of elements in $\phi\upmo(0)$ induced by
the defining equation $\phi$. This identification of two obstruction
theories follows from the canonical isomorphism of 
the C\v{e}ch cohomology and
the Dolbeault cohomology of vector bundles.

The layout of the paper is as follows. In section two, we will recall the analytic
construction of the GW-invariants of symplectic manifolds. We will construct the
Euler class of $[\Phi\mh\bB\to\bE]$ in details using the weakly smooth
approximations constructed in \cite{LT2}. In section three, we will construct a
collection of holomorphic weakly smooth approximations for projective manifolds.
The proof of the Comparisom Theorem will occupy the last section of this paper.   

%

\section{Symplectic construction of GW invariants}


The goal of this section is to review the symplectic construction of 
the GW-invariants of algebraic varieties. We will emphasize on those
parts that are relevant to our proof of the Comparison Theorem.

In this section,we will work mainly with real manifolds and will use the
standard notation in real differential geometry. 

We begin with the symplectic construction of GW-invariants.
Let $X$ be a smooth complex projective variety, and let
$A\in H_2(X,\ZZ)$ and let $g,\, n\in\ZZ$ be fixed once and for all.
We recall the notion of stable $C^l$-maps \cite[Definition 2.1]{LT2}.

\begin{defi}
\lab{1.1}
An $n$-pointed stable map is a
collection $(f;\Sigma;x_1\ldots,x_n)$ satisfying the following property:
First,
$(\Sigma;x_1,\ldots,x_n)$ is an $n$-pointed connected prestable
complex curve with normal crossing
singularity; Secondly, 
$f\mh\Sigma\to X$ is continuous, and the composite $f\circ\pi$ is 
smooth, where $\pi\mh\tilde\Sigma\to\Sigma$ is 
the normalization of $\Sigma$; And thirdly,
if we let $S\sub\Sigma$ be the union of singular locus of $\Sigma$
with its
marked points, then any rational component $R\sub\tilde{\Sigma}$
satisfying $(f\circ \pi)\lsta([R])=0\in H_2(X,\ZZ)$ must
contains at least three points in $\pi\upmo(S)$.
\end{defi}

For convenience, we will abbreviate $(f;\Sigma;x_1,\ldots,x_n)$ to 
$(f;\Sigma;\{x_i\})$. Later, we will use $\cC$ to denote an arbitrary stable
map and use $f_{\cC}$ and $\Sigma_{\cC}$ to denote its corresponding
mapping and domain.
Two stable maps $(f;\Sigma;\{x_i\})$ and $(f\pri;\Sigma\pri;\{x\pri_i\})$
are said to be equivalent if there is an isomorphism 
$\rho\mh\Sigma\to\Sigma\pri$ such that $f\pri\circ\rho=f$
and $x_i\pri=\rho(x_i)$. When 
$(f;\Sigma;\{x_i\})\equiv(f\pri;\Sigma\pri;\{x\pri_i\})$,
such a $\rho$ is called an automorphism of $(f;\Sigma;\{x_i\})$.

We let $\bB$ be the space of equivalence classes $[\cC]$ of
$C^l$-stable maps $\cC$ such that the arithmetic 
genus of $\Sigma_{\cC}$ is $g$ and 
$f_{\cC\ast}([\Sigma])=A\in H_2(X;\ZZ)$.
Note that $\bB$ was denoted by $\bar{\cF}^l_A(X,g,n)$
in~\cite{LT2}.
Over $\bB$ there is a generalized bundle $\bE$ defined as follows.
Let $\cC$ be any stable map and let 
$\tilde f_{\cC}\mh\tilde \Sigma_{\cC}\to X$ be the composite of $f_{\cC}$ with
$\pi\mh\tilde\Sigma_{\cC}\to\Sigma_{\cC}$. We define $\lamzo_{\cC}$
to be the space of all $C^{l-1}$-smooth sections of $(0,1)$-forms of 
$\tilde\Sigma$ with values in ${\tilde f}\sta TX$. Assume $\cC$ and $\cC\pri$
are two equivalent stable maps with $\rho\mh\Sigma_{\cC}
\to\Sigma_{\cC\pri}$ the associated isomorphism, then there is a canonical
isomorphism $\lamzo_{\cC\pri}\cong\lamzo_{\cC}$. We
let $\lamzo_{[\cC]}$ be $\lamzo_{\cC}/\Aut(\cC)$. Then the union 
$$\bE=\bigcup_{[\cC]\in\bB} \lamzo_{[\cC]}
$$
is a fibration over $\bB$ whose fibers are finite quotients of
infinite dimensional linear spaces.
There is a natural section
$$\Phi: \bB\lra\bE
$$
defined as follows. For any stable map $\cC$, we define $\Phi(\cC)$ to
be the image of $\dbar f_{[\cC]}\in\lamzo_{\cC}$ in
$\lamzo_{[\cC]}$. Obviously, for $\cC\sim\cC\pri$ we have
$\Phi(\cC)=\Phi(\cC\pri)$. Thus $\Phi$ descends to a map
$\bB\to\bE$, which we still denote by $\Phi$.

>From now on, we will denote by $\mgna$ the moduli scheme of stable
moprhisms $f\mh C\to X$ with $n$-marked points such that
$C$ is (possibly with nodal singularities) has arithmetic
genus $g$ with $f\lsta([C])=A$.

\begin{lemm}
\lab{1.2}
The vanishing locus of $\Phi$ is canonically homeomorphic to the
underlying topological space of $\mgna$.
\end{lemm}

\begin{proof}
A stable $C^l$-stable map ${\cC}$ in $\bB$
belongs to the vanishing
locus of $\Phi$ if and only if $f_{\cC}$ is holomorphic. 
Since $\Sigma_{\cC}$ is compact,
$\cC$ is the underlying analytic map of a stable morphism. Hence there
is a canonical map $\Phi\upmo(0)\to \mgna\top$,
which is one-to-one and
onto. This proves the lemma.
\end{proof}

To discuss the smoothness of $\Phi$, we need the local uniformizing
charts of $\Phi\mh \bB\to\bE$ near $\Phi\upmo(0)$.  
Let $w\in\bB$ be any point represented by
the stable map $(f_0;\Sigma_0;\{x_i\})$ 
with automorphism group $G_{w}$. 
We pick integers $r_1, r_2>0$ and smooth ample divisors
$H_1,\ldots,H_{r_2}$ with $[H_i]\cdot[A]=r_1$ such that
all $f_0\upmo(H_i)$ are contained in the
smooth locus of $\Sigma_0$ and that
for any $x\in f_0\upmo(H_i)$ we have
\begin{equation}
\image(df_0(x))+T_{f(x)}H_i= T_{f(x)}X.
\lab{eq:1.1}
\end{equation}
Now let $U\sub \bB$ be a sufficiently small neighborhood of
$w\in\bB$ and let
$\tilu$ be the collection of all $(\cC; z_{n+1},\ldots,z_{n+r_1r_2})$
such that $\cC\in U$ and the $z_i$'s
is a collection of smooth points of $\Sigma_{\cC}$ such that
for each $1\leq j\leq r_2$ the subcollection
$(z_{n+(j-1)r_1+1},\ldots,z_{n+jr_1})$ contains distinct points
and is exactly $f_{\cC}\upmo(H_j)$.
Note that we do not require $(z_{n+1},\ldots,z_{n+r_1r_2})$ to be distinct.
\footnote{In case $X$ is a symplectic manifold, then we should use locally
closed real codimension 2 submanifold instead of $H_i$, as
did in \cite{LT2}. 
Here we use this construction of uniformizing charts
because it is compatible to the construction of atlas of the
stack $\mgna$ in algebraic geometry.}
Let $\pitilu\mh\tilu\to U$ be the projection 
that sends $(\cC;z_{n+1},\ldots,z_{n+r_1r_2})$ to $\cC$.
Clearly, $G_w$ acts on $\pi_U\upmo(w)$ canonically by permuting
their $(n+r_1r_2)$-marked points. Namely, for any
$\sigma\in G_w$ and $\cC\in\pi_U\upmo(w)$ with 
marked points $z_1,\ldots,z_{n+r_1r_2}$, $\sigma(\cC)$ is the
same map with the marked points $\sigma(z_1),\ldots,\sigma(z_{n+r_1r_2})$.
In particular, we can view $G_w$ as a subgroup 
of the permutation group $S_{n+r_1r_2}$. 
Hence $G_w$ acts on $\tilde U$ by permuting the marked
points of $\cC\in\tilde U$ according to the inclusion $G_w\sub S_{n+r_1r_2}$.
Note that if $H_i$'s are
in general position then elements in $\tilu$ has no automorphisms
and have distinct marked points.
Let $G\ltilu=G_w$. Since fibers of $\pi_U$ are invariant
under $G\ltilu$, $\pi_U$ induces a map $\tilu/G\ltilu\to U$,
which is obviously a covering\!
\footnote{In this paper we call $p\mh A\to B$ a covering if $p$ is a
covering projection~\cite{Sp}
and $\#(p\upmo(x))$ is independent of $x\in B$.
We call $p\mh A\to B$ a local covering if $p(A)$ is open in $B$ and
$p\mh A\to  p(A)$ is a covering.} if $U$ is
sufficiently small. Further, if we let
$$\bE\ltilu=\bigcup_{\cC\in\tilu} \lamzo_{\cC}
$$
and let $\Phi\ltilu\mh\tilu\to\bE\ltilu$ be the section that sends
$\cC$ to $\dbar f_{\cC}$, then $\Phi\ltilu$ is 
$G\ltilu$-equivariant and $\Phi|_U\mh U\to \bE|_U$ is the descent of
$\Phi\ltilu/G\ltilu\mh\tilu/G\ltilu\to\bE\ltilu/G\ltilu$.
Note that fibers of $\bE\ltilu$ over $\tilu$ are linear spaces. Following
the convention, we will call 
$\Lambda=(\tilu,\etilu,\phitilu,\gtilu)$ a uniformizing chart of
$(\bB,\bE,\Phi)$ over $U$. Let $V\sub U$ be an open subset
and let $\tilde V=\pi_U\upmo(V)$, let $G\ltilv=G\ltilu$,
let $\bE\ltilv=\bE\ltilu|_{\ltilv}$ and let $\Phi\ltilv=\Phi\ltilu|\ltilv$.
We will call $\Lambda\pri=(\tilde V,\bE\ltilv,\Phi\ltilv, G\ltilv)$ a
uniformizing chart of $(\bB,\bE,\Phi)$ that is 
the {\it restriction} of the original chart to $V$, and denoted by $\Lambda|_V$.
We can also construct uniformizing charts by pull back. Let
$G\ltilv$ be a finite group acting effectively on
a topological space $\tilde V$, let $G\ltilv\to G\ltilu$
be a homomorphism and $\varphi\mh \tilde V\to\tilde U$ be a
$G\ltilv$-equivariant map so that $\tilde V/G\ltilv\to \tilde U/G\ltilu$
is a local covering map. Then we set $\bE\ltilv=\varphi\sta\bE\ltilu$ and
$\Phi\ltilv=\varphi\sta\Phi\ltilu$. The data
$\Lambda\pri=(\tilde V,\bE\ltilv,\Phi\ltilv, G\ltilv)$
is also a uniformizing chart. We will call $\Lambda\pri$ the
{\sl pull back} of $\Lambda$, and denoted by $\varphi\sta\Lambda$.
In the following, we will denote the collection of all uniformizing charts of
$(\bB,\bE,\Phi)$ by $\C$.

The collection $\C$ has the following compatibility
property. Let 
$$\Lambda_i=(\tilu_i,\bE_{\tilu_i},\Phi_{\tilu_i}, G_{\tilu_i}),
$$
where $i=1,\ldots,k$, be a collection of uniformizing charts in $\C$
over $U_i\sub\bB$ respectively. Let $p\in\cap_{i=1}^k U_i$
be any point. Then there is a uniformizing chart
$\Lambda=(\tilv,\bE_{\tilv},\Phi_{\tilv},G_{\tilv})$ over $V\sub\cap^k U_i$ 
with $p\in V$
such that there are homomorphisms $G_{\tilv}\to G_{\tilu_i}$
and equivariant local covering maps $\varphi_i\mh\tilv\to\pi_{U_i}\upmo(V)
\sub\tilu_i$ compatible with $\tilv\to V$ and 
$\pi_{U_i}\upmo(V)\to V\sub U_i$, such that
$\varphi_i\sta(\bE_{\tilu_i},\Phi_{\tilu_i})
\cong (\bE_{\tilv},\Phi_{\tilv})$.
In this case, we say $\Lambda$ is  
{\it finer than} 
$\Lambda_i|_{V}$.

The main difficulty in constructing the GW invariants in this setting is that
the smoothness of $(\tilu,\etilu,\phitilu)$ is unclear when $U$ 
contains maps whose domains are singular.
To overcome this difficulty, the authors introduced the notion of generalized
Fredholm bundles in~\cite{LT2}. The main result of~\cite{LT2} 
is the following theorems, 
which enable them to construct the GW invariants for all symplectic manifolds.

\begin{theo}
\lab{1.3}
The data $[\Phi\mh \bB\to\bE]$ is a generalized oriented Fredholm
V-bundle of relative index $2r\lexp$, where $r\lexp=
c_1(X)\cdot A+n+(n-3)(1-g)$ is half of the virtual (real) dimension of
$\Phi\upmo(0)$.
\end{theo}

\begin{theo} 
\lab{1.4}
For any generalized oriented Fredholm V-bundle
$[ \Phi\mh \bB\to\bE]$ of relative index $r$, we can assign
to it an Euler class $e([\Phi\mh \bB\to\bE])$ in $H_r(\bB;\QQ)$ 
that satisfies all the expected properties of the Euler classes.
\end{theo}

As explained in the introduction, the pairing of the Euler class of
$[\Phi\mh \bB\to\bE]$ with the tautological topological class 
will give rise to the symplectic version of the GW invariants of $X$. 
Further, the Comparison Theorem we set out to prove amounts to compare
this Euler class with the image of the virtual moduli cycle $[\mgna]\vir$
in $H_r(\bB;\QQ)$ via the inclusion $\mgna\top\sub\bB$.
In the remainder part of this section, we will list all properties of 
$ [\Phi\mh \bB\to\bE] $ that are relevant to the construction of its
Euler class. This list is essentially equivalent to saying
that $ [\Phi\mh \bB\to\bE] $ is a generalized oriented Fredholm V-bundle. 
After that, we will construct the Euler class of $ [\Phi\mh \bB\to\bE]  $ 
in details.


We begin with the notion of weakly smooth structure. A local smooth 
approximation of $[\Phi\mh \bB\to\bE]$ over $U\sub \bB$
is a pair $(\Lambda,V)$, where $\Lambda=(\tilu,\etilu,\phitilu,G\ltilu)$
is a uniformizing chart over $U$ and $V$ is a finite equi-rank
$G\ltilu$-vector bundle over $\tilu$
that is a $\gtilu$-equivariant subbundle of $\etilu$ such that
$R_V\eqdef \phitilu\upmo(V)\sub\tilu$ is an equi-dimensional
smooth manifold, $V|_{R_V}$ is a smooth vector bundle and the lifting
$\phi_V\mh R_V\to V|\lrv$ of $\phitilu|\lrv$ is a smooth section.
An orientation of $(\Lambda,V)$ is a $G_{\tilv}$-invariant
orientation of the real line bundle $\wedge\top(TR_V)\otimes\wedge\top
(V|_{R_V})\upmo$ over $R_V$.
We call $\rank V-\dim R_V$ the
index of $(\Lambda,V)$ (We remind that all ranks and dimensions
in this section are over reals). Now
assume that $(\Lambda\pri,V\pri)$ is another weakly smooth structure of
identical index over $W\sub\bB$. We say that $(\Lambda\pri,V\pri)$
is finer than $(\Lambda,V)$ if the following holds.
First, the restriction $\Lambda\pri|_{W\cap U}$ is finer than
$\Lambda|_{W\cap U}$; Secondly, if we let $\varphi\mh
\pi_W\upmo(W\cap U)\to\pi_U\upmo(W\cap U)$
be the covering map then
$\varphi\sta V\sub \varphi\sta\etilu\equiv
\bE\ltilw|_{\pi_W\upmo(W\cap U)}$
is a subbundle of $V\pri|_{\pi_W\upmo(W\cap U)}$;
Thirdly, for any $w\in\tilw$ the homomorphism $
T_wR_{V\pri}\to \bl V\pri/\varphi\sta V\br |_w$ induced by
$d\phi_{V\pri}(w)\mh T_wR_V\to V\pri|_w$ is surjective, and
the map $\phi_{V\pri}\upmo(\varphi\sta V)\to R_V$ 
induced by $\varphi$ is
a local diffeomorphism between smooth manifolds.
Note that the last condition implies that if we identify 
$T_{\varphi(w)}R_V$ with 
$T_w\phi_{V\pri}\upmo(\varphi\sta V)
\sub T_w R_{V\pri}$, then the induced homomorphism
\begin{equation}
\lab{eq:1.2}
T_wR_{V\pri}
/T_{\varphi(w)}R_V\lra \bl V\pri/\varphi\sta V)|_w
\end{equation}
is an isomorphism. In case both $(\Lambda,V)$ and $(\Lambda\pri,V\pri)$
are oriented, then we require that the orientation of 
$(\Lambda,V)$ coincides with that of $(\Lambda\pri,V\pri)$
based on the isomorphism
\begin{equation}
\lab{eq:1.3}
\wedge\top(T_w R_{V\pri})\otimes\wedge\top(V\pri|_w)\upmo\cong
\wedge\top(T_{\varphi(w)}R_V)\otimes\wedge\top(V|_{\varphi(w)})\upmo
\end{equation}
induced by~\eqref{eq:1.2}.

Now let $\A=\{(\Lambda_i,V_i)\}_{i\in\cK}$ be a collection of oriented 
smooth approximations of $(\bB,\bE,\Phi)$. In the following, we will denote
by $U_i$ the open subsets of $\bB$ such that $\wedge_i$ is a 
smooth chart over $U_i$. We say $\A$ covers $\phi\upmo(0)$ if 
$\phi\upmo(0)$ is contained in the union of the images of $U_i$ in $\bB$.

\begin{defi}\lab{1.5}
An index $r$ oriented weakly smooth structure of $(\bB,\bE,\Phi)$ is a 
collection $\A=\{(\Lambda_i,V_i)\}_{i\in\cK}$ of index $r$ oriented smooth
approximations such that $\A$ covers $\Phi\upmo(0)$ and that
for any $(\Lambda_i,V_i)$ and $(\Lambda_j,V_j)$ in $\A$ with 
$p\in U_i\cap U_j$, there is a $(\Lambda_k,V_k)\in\A$ such that 
$p\in U_k$ and $(\Lambda_k,V_k)$ is finer than $(\Lambda_i,V_i)$ and
$(\Lambda_j,V_j)$.
\end{defi}

Let $\A\pri$ be another index $r$ oriented weakly smooth structure of
$(\bB,\bE,\Phi)$. We say $\A\pri$ is finer than $\A$ if for any
$(\Lambda,V)\in\A$ over $U\sub\bB$ and $p\in U\cap\Phi\upmo(0)$,
there is a $(\Lambda\pri,V\pri)\in\A\pri$ over $U\pri$ such that
$p\in U\pri$ and $(\Lambda\pri,V\pri)$ is finer than
$(\Lambda,V)$.
We say that two weakly smooth structures $\A_1$ and $\A_2$ 
are equivalent if there is a third weakly smooth structure that is
finer than both $\A_1$ and $\A_2$.

\begin{prop}[\cite{LT2}]\lab{1.6}
The tuple $(\bB,\bE,\Phi)$ constructed at the beginning of this section
admits a canonical oriented weakly smooth structure of index $2r\lexp$.
\end{prop}

We remark that the construction of such a weakly smooth structure is
the core of the analytic part of [LT2]. 

In the following, we will use the weakly smooth structure
of $[\Phi\mh\bB\to\bE]$ to construct its Euler class.
The idea of the construction is as
follows. Given a local smooth approximation $(\Lambda,V)$ over
$U\sub\bB$, we obtain a smooth manifold $R_V$, a vector bundle $V|\lrv$ 
and a smooth section $\phi_V\mh R_V\to V|\lrv$. Following the 
topological construction of the Euler classes, we shall perturb 
$\phi_V$ to a new section $\tilde\phi_V\mh R_V\to V|\lrv$ so that
$\tilde\phi_V$ is transversal to the zero section of $V|\lrv$. Here
by a section transversal to the zero section, we mean that the
graph of this section is transversal to the zero section in the
total space of the vector bundle. Hence the 
Euler class will be the cycle represented by ${\tilde\phi_V}\upmo(0)$
near $U$. Since the weakly smooth structure of $\Phi\mh \bB\to\bE$ 
is given by a collection of compatible by not necessary matching local
smooth approximations, we need to work out this perturbation scheme 
with special care so that $\{\tilde\phi_V\upmo(0)\}$ patch together
to form a well-defined cycle.

Let $\A=\{(\Lambda\lalp,V\lalp)\}_{\alpha\in\cK}$ be the weakly smooth 
structure provided by Proposition~\ref{1.6}. For convenience, for any
$\alpha\in\cK$ we will denote the corresponding uniformizing chart 
$\Lambda\lalp$ by $(\tilu\lalp,\etil\lalp,\Phitil\lalp,G\lalp)$ and
will denote its descent by $(U\lalp,E\lalp,\Phi\lalp)$.
Accordingly, we will denote the projection $\pi_{U\lalp}\mh\util\lalp\to U\lalp$ 
by $\pi\lalp$, denote ${\tilde\Phi\lalp}\upmo(V\lalp)$ by $R\lalp$, denote 
$V\lalp|_{R\lalp}$ by $W\lalp$ and denote the lifting of $\tilde\Phi\lalp|_{R\lalp}
\mh R\lalp\to\tile\lalp|_{R\lalp}$ by $\phi\lalp\mh R\lalp\to W\lalp$.
Without loss of generality, we can assume that for any approximation
$(\Lambda\lalp,V\lalp)\in\A$ over $U\lalp$ and any
$U\pri\sub U\lalp$, the restriction $(\Lambda\lalp,V\lalp)|_{U\pri}$
is also a member in $\A$.
In the following, we call $S\sub R\lalp$ {\it symmetric} if 
$S=\pi\lalp\upmo(\pi\lalp(S))$.

Next, we pick a covering data for $\Phi\upmo(0)\sub\bB$
provided by the following covering lemma.

\begin{lemm}[\cite{LT2}]\lab{1.20}
There is a finite collection $\cL\sub\cK$ and a total ordering of
$\cL$ of which the following holds.
the set $\Phi\upmo(0)$ is contained in the union of 
$\{R\lalp\}_{\alpha\in\cL}$ and
for any $\alpha$ and $\beta\in\cL$ such that $\alpha<\beta$
then approximation $(\Lambda\lbe,V\lbe)$ is finer than the
approximation $(\Lambda\lalp,V\lalp)$.
\end{lemm}

\begin{proof}
The lemma is part of Proposition 2.2 in \cite{LT2}. It is proved there
by using the stratified structures of $(\bB,\bE,\Phi)$.
Here we will give a direct proof of this by using the definition of
smooth approximations, when $\Phi\upmo(0)$ is triangulable,
which is true when $X$ is projective.
Let $k$ be the real dimension of 
$\Phi\upmo(0)$. To prove the lemma, we will show that
there are $k+1$ subsets $\cL_k,\ldots,\cL_0\sub\cK$ and that
for each $\alpha\in\cup_{i=0}^k\cL_i$ 
there is an open symmetric subset $U\lalp\pri\subb U\lalp$
such that $R\lalp\pri=R\lalp\cap\pi\lalp\upmo(U\lalp\pri)\subb R\lalp$
of which the following holds: first, for each $i\leq k$ the set 
$Z_i=\Phi\upmo(0)-\cup_{j\geq i}\cup_{\alpha\in\cL_j} U\lalp\pri$
is a triangulable space whose dimension is at most
$i-1$, and secondly, for any pair of distinct 
$(\alpha,\beta)\in\cL_i\times\cL_j$ 
with $i\leq j$, the restriction 
$(\Lambda\lalp,V\lalp)|_{U\lalp\pri\cap U\lbe\pri}$
is finer than $(\Lambda\lbe,V\lbe)|_{U\lalp\pri\cap U\lbe\pri}$.
We will construct $\cL_i$ inductively, starting from $\cL_k$. We first pick
a finite $\cL_k\sub\cK$ so that $\cup_{\alpha\in\cL_k} U\lalp
\supset \Phi\upmo(0)$. This is possible since 
$\Phi\upmo(0)$ is compact. Since it is also triangulable, we can
find a symmetric $U\lalp\pri\sub U\lalp$ for each $\alpha\in\cL_k$
so that $\{U\pri\lalp\}_{\alpha\in\cL_k}$ is disjoint,
$R\lalp\pri\subb R\lalp$ and
$Z_k$
is trangulable with dimension at most $k-1$.
Now we assume that we have found $\cL_k,\ldots,\cL_i$ as desired.
Then for each $x\in Z_i$ we can find a neighborhood $O$ of
$x\in\bB$ such that for any $\alpha\in\cup_{j\geq i}\cL_j$ either
$x\in U\lalp$ or $O\cap U\lalp\pri=\emptyset$. Let $\cI_x$ be those
$\alpha$ in $\cup_{j\geq i}\cL\lalp$ such that $x\in U\lalp$. Then by the
property of $\A$ there is a $\beta\in\cK$ so that $(\Lambda\lbe,V\lbe)$ is
finer than $(\Lambda\lalp,V\lalp)|_{U\lalp\pri}$ for all $\alpha\in\cI_x$.
Without loss of generality, we can assume that $U\lbe\sub O$. Then
$(\Lambda\lbe,V\lbe)$ is finer than $(\Lambda\lalp,V\lalp)$ for
all $\alpha\in\cup_{j\geq i}\cL_j$.
Since $Z_i$ is compact, we can cover it by finitely many
such $(\Lambda\lbe,V\lbe)$'s, say indexed by
$\cL_{i-1}\sub\cK$.
On the other hand, since 
$Z_i$ is triangulable with dimension at most $i-1$, we can find symmetric
$U\lalp\pri\sub  U\lalp$ for each $\alpha\in\cL_{i-1}$ so that 
$R\pri\lalp\subb R\lalp$ for $\alpha\in\cL_{i-1}$ and
$Z_i-\cup_{\alpha\in\cL_{i-1}} U\lalp\pri$ is trianglable
with dimension at most $i-2$. This way, we can find the set $\cL_k,\ldots,
\cL_0$ as desired. In the end, we simply put $\cL=\cup_{i=0}^k\cL_i$.
We give it a total ordering so that whenever $\alpha\in\cL_i$, $i\geq j$
and $\beta\in\cL_j$ then $\alpha\leq\beta$. This proves the Lemma.
\end{proof}

We now fix such a collection $\cL$ once and for all. Since
$\cL$ is totally ordered, in the following we will replace the
index by integers that range from 1 to $\#(\cL)$ and use
$k$ to denote an arbitrary member of $\cL$.
We first build the comparison data into the collection
$\{R_k\}_{k\in\cL}$ and $\{W_k\}_{k\in\cL}$. To distinguish the
projection $\pi_k\mh\util_k\to U_k$ from the composite
$\util_k\to U_k\to\bB$, we will denote the later by $\iota_k$.
For any pair $k\geq l$, we set $R\lkl=\iota_k\upmo(\iota_l(R_l))$.
Then there is a canonical map and a 
canonical vector bundle inclusion
\begin{equation}
\lab{eq:1.4}
f\ulk\mh R\lkl\to R_l
\and
(f\ulk)\sta(W_l)\mapright{\sub} W_k|_{R\lkl},
\end{equation}
that is part of the data making $(\Lambda_k,V_k)$
finer than $(\Lambda_l,V_l)$. 
Note that $R_{k,l}\sub R_k$ is a
locally closed submanifold,
$f^l_k(R\lkl)$ is open in $R_l$ and $f\ulk\mh R\lkl\to f\ulk(R\lkl)$
is a covering map. Because of the compatibility condition, for any 
$k>l>m$ if $R\lkl\cap R_{k,m}\ne\emptyset$ then
$f\ulk(R\lkl\cap R_{k,m})\sub R_{l,m}$ and 
\begin{equation}
f^m_l\circ f^l_k=f^m_k: 
R\lkl\cap R_{k,m}\lra R_m.
\lab{eq:1.5}
\end{equation}
Further, restricting to $R_{k,l}\cap R_{k,m}$,
the pull backs
\begin{equation}
\lab{eq:1.6}
(f_k^m)\sta(W_m)|_{R_{k,l}\cap R_{k,m}}
=(f\ulk)\sta(f^m_l)\sta(W_m)|_{R_{k,l}\cap R_{k,m}}
\sub W_k|_{R_{k,l}\cap R_{k,m}}.
\end{equation}
In the following, we will use $\R$ to denote the collection
of data $\{(R\lkl,f\ulk)\}$ and use $\W$ to denote
the data $\{(W_k,(f\ulk)\sta)\}$. We will call the pair $(\R,\W)$
a good atlas
of the weakly smooth structure $\A$ of $[\Phi\mh\bB\to\bE]$.
For technical reason, later we need to shrink each $R_k$ slightly.
More precisely, let $\{ S_k\}_{k\in\cL}$ be a collection of 
symmetric open subsets $S_k\subb R_k$ such that $\{ S_k\}$ 
still covers $\Phi\upmo(0)$. We then let $S\lkl=(f^l_k)\upmo(S_l)\cap
S_k$, let $W_k\pri=W_k|_{S_k}$ and let $g^l_k$ and $(g^l_k)\sta$ be the
restriction to $S\lkl$ of $f^l_k$ and $(f^l_k)\sta$ respectively. 
Then $(\fS,\W\pri)$, where $\fS=\{(S\lkl,g^l_k)\}$ and $\W\pri=
\{(W_k\pri,(g^l_k)\sta)\}$, is also a good atlas of $[\Phi\mh \bB\to\bE]$.
We call it a precompact sub-atlas of $(\R,\W)$, and denote it 
in short by $\fS\subb\R$.

%


To describe the collection $\{\phi_k\}$, we need to introduce the notion of
regular extension. Let $M$ be a manifold and $M_0\sub M$
be a locally closed submanifold. Let $V\to M$ be a smooth
vector bundle and $V_0\to M_0$ a subbundle of $V|_{M_0}$.
We assume that both $(M,V)$ and $(M_0,V_0)$ are 
oriented. We say that a section $h\mh M\to V$ is a smooth extension of 
$h_0\mh M_0\to V_0$ if both $h_0$ and $h$ are smooth and if
the induced section $M_0\mapto{h_0} V_0\to V|_{M_0}$
is identical to the restriction $h|_{X_0}\mh X_0\to V_0$. 
We say $h$ is a regular extension of $h_0$ if in addition to $h$
being a smooth extension of $h_0$ we have that for any $x\in X_0$ the 
homomorphism
\begin{equation}
dh(x): T_xM/T_xM_0\lra (V/V_0)|_x
\lab{eq:1.7}
\end{equation}
is an isomorphism and the orientation of $(M,V)$ and $(M_0,V_0)$ are
compatible over $M_0$ based on the isomorphism~\eqref{eq:1.7}.

\begin{defi}
A collection $\{h_k\}_{k\in\cL}$ is called a smooth section of $\W$ if
$h_k$ is a smooth section of $W_k$ for each $k\in\cL$ and $h_k$
is a smooth extension of $h_l$ for any pair $k\geq l$ in $\cL$.
If in addition that $h_k$ is a regular extension of
$h_l$ for all $k\geq l$, then we call $\{h_k\}$ a regular section of
$\W$.
\end{defi}

In the following, we will use $\h\mh\R\to\W$ to denote a 
smooth section with $\h$ understood to be $\{h_k\}_{k\in\cL}$.
We set $\h\upmo(0)$ to be the collection $\{h_k\upmo(0)\}$
and set $\iota(\h\upmo(0))$ to be the union of $\iota_k(h_k\upmo(0))$
in $\bB$.
We say $\h\upmo(0)$ is proper if 
$\iota(\h\upmo(0))$ is compact. 
Without loss of generality, we can assume that $\dim R_k>0$ for all
$k\in\cL$.
We say that $\h$ is transversal to the zero section $\mathbf 0\mh \R\to\W$
if $\h$ is a regular section and
if for any $k\in\cL$ the graph $\Gamma_{h_k}$ of $h_k$ is 
transversal to the 0 section of $W_k$ in the total space of $W_k$.

\begin{lemm}
\lab{1.9}
Let the notation be as before. Then
$\h\upmo(0)$ is proper if and only if there is a symmetric open
subsets $R_k\pri\subb R_k$ for each $k\in\cL$ such that 
$\cup_{k\in\cL}\iota_k(h_k\upmo(0))\sub
\cup_{k\in\cL}\iota_k(R\pri_k)$
and such that for each $k\in\cL$,
\begin{equation}
\lab{eq:1.8}
h_k\upmo(0)\cap (R_k- R_k\pri)
\sub\bl\bigcup_{l<k}(f\ulk)\upmo(R_l\pri)\br
\cup\bl \bigcup_{l>k}f\ulk(R\lkl\pri)\br.
\end{equation}
\end{lemm}

\begin{proof}
We first assume that $Z=\cup_{k\in\cL} \iota_k(h_k\upmo(0))$ is compact. 
Then since $\{R_k\}_{k\in\cL}$ covers $Z$ and since $\dim R_k>0$,
for each $k\in\cL$ we can find symmetric $R_k\pri\subb R_k$  so that 
$\{R_k\pri\}_{k\in\cL}$
still covers $Z$. Obviously, this implies~\eqref{eq:1.8}.
Conversely, if we have found $R_k\pri\sub R_k$ 
as stated in the lemma, then  $\{\clo(\iota_k(R_k\pri))\cap Z\}$ will
cover $Z$, where $\clo(A)$ is the closure of $A$. 
Since $\clo(\iota_k(R_k\pri))$ are compact and since
$Z\cap \clo(\iota_k(R_k\pri))$ is closed in $\clo(\iota_k(R_k\pri))$, $Z$ is
compact as well. This proves the lemma.
\end{proof}

\begin{lemm}
\lab{1.10}
Let $\bphi\mh \R\to\W$ be the collection $\{\phi_k\}$ 
induced by $\{\tilde\Phi_k\}_{k\in\cL}$. Then $\bphi$ is
a regular section with proper vanishing locus.
\end{lemm}

\begin{proof}
This is equivalent to the fact that $[\Phi\mh\bB\to\bE]$ is a 
weakly Fredholm V-bundle, which was introduced and proved in
\cite{LT2}.
\end{proof}

Now let $\h\mh\R\to\W$ be a regular section such that $\h$ is 
transversal to the zero section and $\h\upmo(0)$ is proper.
We claim that the data $\{h_k\upmo(0)\}$ descends to an oriented current
in $\bB$
with rational coefficients supported on a stratified set whose
boundary is empty. In particular, it defines a singular homology 
class in $H\lsta(\bB,\QQ)$.

Recall that for each $k\in\cL$ the associated group $G_k$ acts on
$R_k$ such that $R_k/G_k$ is a covering of $\iota_k(R_k)$. 
We let $m_k$ be the product of the order of $G_k$ 
with the number of the sheets
of the covering $R_k/G_k\to\iota_k(R_k)$. Note that
then the covering $R\lkl\to f\ulk(R\lkl)$ is an $m_k/m_l$-fold
covering. Because $h_k$ is a regular extension of 
$(f\ulk)\sta(h_l)$, $(f\ulk)\sta(h_l)\upmo(0)$
is an open submanifold of $h_k\upmo(0)$ with identical orientations.
Hence $\iota_k(h\upmo_k(0))$ and $\iota_l(h_l\upmo(0))$
patch together to form a stratified subset, and consequently
the collection $\{\iota_k(h_k\upmo(0))\}_{k\in\cL}$
patch together to form a stratified subset, say $Z$, in $\bB$.
Now we assign multiplicities to open strata of $Z$. 
Let $O_k=\iota_k(h_k\upmo(0))$. Since $O_k\sub Z$ is
an open subset, we can assign multiplicities to $O_k$
so that as oriented current $[O_k]=\iota\lsta(\frac{1}{m_k}[
h_k\upmo(0)])$, where $[h_k\upmo(0)]$ is the current of 
the oriented manifold $h_k\upmo(0)$ with multiplicity one.
Here the orientation of $h_k\upmo(0)$ is the one induced
by the orientation of $(R_k,W_k)$.
Using the fact that $R\lkl\to f\ulk(R\lkl)$ is a covering
with $m_k/m_l$ sheets, the assignments of the multiplicities 
of $O_k$ and $O_l$ over $\iota_k(R_k)\cap\iota_l(R_l)$ coincide.
Therefore $Z$ is an oriented stratified set of pure dimension 
with rational multiplicities. We let $[Z]$ be the corresponding current.
It remains to check that $\partial[Z]=0$ as current.
Clearly, $\partial[Z\cap O_k]\sub\clo(O_k)-O_k$. Since
$\{O_k\cap Z\}$ is an open covering of $Z$, $\partial[Z]=0$ if
$Z$ is compact. But this is what we have assumed in the first place.
Later, we will denote the so constructed cycle by
$$[\h\upmo(0)]\in H\lsta(\bB,\QQ).
$$

In the remainder of this section, we will perturb the section
$\bphi\mh\R\to\W$ to a new section so that it is transversal to the zero section
and so that its vanishing locus is compact. 
The current defined by the vanishing locus of the perturbed section
will define the Euler class of $[\Phi\mh\bB\to\bE]$.
 
We begin with a collection $\cS=\{S_l\}_{l\in\cL}$
of symmetric open $S_k\subb R_k$ such that $\{\iota_l(S_l)\}$ 
cover $\iota(\bphi\upmo(0))$.
For technical reason, we assume that for each $k\in\cL$ the
boundary $\partial S_k$, which is defined to be $\clo(S_k)-S_k$
in $R_k$, is a 
smooth manifold of dimension $\dim S_k-1$. By slightly altering
$S_k$ if necessarily, we can and do assume that $\partial S_k$ is 
transversal to $R\lkl$ along $\partial
S_k\cap (f\ulk)\upmo(\clo(S_l))$ for all $l<k$. (We will call
such $\cS$ satisfying the transversality condition on
its boundary.) Following the convention, 
we set $S\lkl=(f\ulk)\upmo(S_l)\cap S_k$. We now construct a collection of 
(closed) tubular neighborhoods of $S\lkl$ in $R_k$. We fix the
index $k$ and consider the closed submanifold 
(with boundary) $\Sigma_l:=\clo(S\lkl)\sub R_k$. Because of the
transversality condition on $\partial S_l$ and on $\partial S_k$,
we can find a $D^h$-bundle $p_l\mh T_l\to \Sigma_l$,
where $D^h$ is the closed unit ball in $\RR^h$ and
$h=\dim R_k-\dim R\lkl$, and a smooth embedding
$\eta_l\mh T_l\to R_k$ of which the following two conditions
holds. First, the restriction
of $\eta_l$ to the zero section $\Sigma_l\sub T_l$ is the original
embedding $\Sigma_l\sub R_k$, and secondly
\begin{equation}
\lab{eq:1.9}
\eta_l(p_l\upmo(\Sigma_l\cap \partial S_k ))\sub \partial S_k 
\and
\eta_l(p_l\upmo(S_{k,l} ))\sub S_k.
\end{equation}
For any $0<\eps<1$, we let $T_l\ueps\sub T_l$ be the closed
$\eps$-ball subbundle of $T_l$. By abuse of notation, in
the following we will not distinguish $T\ueps_l$
from its image $\eta_l(T_l\ueps)$
in $R_k$. We will call $T\ueps_l$ the 
$\eps$-tubular neighborhood of $\Sigma_l$ in $R_k$. 
One property we will use later is that if $R\lkl\cap R_{k,l\pri}\ne\emptyset$
for $l\pri<l<k$, then $R_{k,l}\cap R_{k,l\pri}$ is an open
subset of $R_{k,l\pri}$, and hence for $0<\eps\ll 1$ we have
$\Sigma_{l\pri}\cap T\ueps_l\sub \Sigma_{l\pri}\cap \Sigma_l$.

Now consider $\Sigma_l\sub R_k$. Since $T_l$ is a disk bundle
over $\Sigma_l$, it follows that we can extend
the subbundle $(f^l_k)\sta(W_l)|_{\Sigma_l}\sub W_k|_{\Sigma_l}$
to a smooth subbundle of $W_k|_{T_l}$, denoted by
$F_l\sub W_k|_{T_l}$.
We then fix an isomorphism and the inclusion
\begin{equation}
p_l\sta\bl (f_k^l)\sta(W_l)|_{\Sigma_l}\br
\cong F_l.
\lab{eq:1.10}
\end{equation}
In this way, we can  extend any section $\zeta$ of 
$(f_k^l)\sta(W_l)|_{\Sigma_l}$ to a section of
$W_k|_{T_l}$ as follows. We first let $\zeta\pri\mh T_l\to F_l$
be the obvious extension 
using the isomorphism \eqref{1.10}. We then let 
$\zeta\lex\mh T_l\to W_k|_{T_l}$
be the induced section using the inclusion 
$F_l\sub W_k|_{T_l}$.
We will call $\zeta\lex$ the standard extension of $\zeta$ to $T_l$.
We
fix a Riemannian metric on $R_k$ and a metric on $W_k$. For any
section $\zeta$ as before, we say $\zeta$ is sufficiently small if its
$C^2$-norm is sufficiently small.
We now state a simple but important observation. 

\begin{lemm}
\lab{1.11}
Let the notation be as before. Then
there is an $\eps>0$ such that for any section 
$g\mh T_l\to F_l\sub W_k|_{T_l}$
such that $\parallel\! g\!\parallel_{C^2}<\eps$, the section $h_k|_{T_l}+g$
is non-zero over $T\ueps_l-\Sigma_l$.
\end{lemm}

\begin{proof}
This follows immediately from the fact that 
$\Sigma_l$ is compact and that for any $x\in R_{k,l}$ the
differential
$$dh_k: T_xR_k/T_x R_{k,l}\lra \bl W_k / (f^l_k)\sta(W_l)\br|_x
$$
is an isomorphism.
\end{proof}

We now state and prove the main proposition of this section.

\begin{prop}
\lab{1.12}
Let $\h\mh \R\to\W$ be a regular section with $\h\upmo(0)$
proper, let $\R\pri\subb\R$ be a good sub-atlas and
let $\h\pri$ be the the restriction of $\h$ to $\R\pri$.
We assume that the vanishing locus of $\h\pri$ is still proper.
Then there is a smooth family of
regular sections $\g(t)\mh \R\pri\to\W\pri$,
where $\W\pri$ be the restriction of $\W$ to $\R\pri$,
parameterized by $t\in [0,1]$ such that
\begin{equation}
\bigcup_{t\in[0,1]}\iota\bl\g(t)\upmo(0)\br\times\{t\}
\sub \bB\times [0,1]
\lab{eq:1.11}
\end{equation}
is compact, that $\g(0)=\h\pri$
and that $\g(1)$ is transversal to the zero section of $\W\pri$.
\end{prop}

\begin{proof}
We will construct the perturbation over $R_1\pri$ and then
successively extends it to the remainder of $\{R\pri_k\}$.
We first fix 
a collection of symmetric open subsets 
$\{S_k\}_{k\in\cL}$ such that $R_k\pri\subb S_k\subb R_k$ and
that $S_k$ satisfies the transversality condition on its boundary. 
Let $k$ be any positive
integer no bigger than $\#(\cL)+1$. The induction hypothesis $\cH_k$ 
states that for each integer $l<k$ we have constructed a symmetric 
open $S_l\pri$ satisfying 
$R_l\pri\subb S_l\pri\subb S_l$ and a smooth family of small 
enough sections $e_l(t)\mh R_l\to W_l$ such that 
$e_l(0)\equiv0$ of which the following holds.
First, let $\bh_l(t)=h_l+e_l(t)$, then for any $l<m<k$ the section
$\bh_{m}(t)|_{S_m\pri}$ is a regular extension of
$(f^l_{m})\sta(\bh_l(t))|_{S_{m,l}\pri}$;
Secondly, for any $l<k$, the section $\bh_l(1)$ is 
transversal to the zero section of $W_l$
over $S_l\pri$, and finally, for any $l<k$ and $t\in[0,1]$,
\begin{equation}
\bh_l(t)\upmo(0)\cap\bl S_l\pri-R_l\pri\br\sub
\bl\bigcup_{i\leq l}(f_l^i)\upmo(R_i\pri)\br\bigcup
\bl\bigcup_{m\geq l}f_m^l(R_{m,l}\pri)\br.
\lab{eq:1.12}
\end{equation}

Clearly, the condition $\cH_1$ is automatically satisfied.
Now assume that we have found $\{ S_l\pri\}_{l<k}$ and
$\{e_l\}_{l<k}$ required by the condition $\cH_k$. We will demonstrate how 
to find $e_k$ and a new sequence of open subsets
$\{S_l\pri\}_{l\leq k}$ so that the condition $\cH_{k+1}$ will hold for
$\{e_l\}_{l\leq k}$ and $\{ S\pri_l\}_{l\leq k}$.

We continue to use the notation developed earlier. In particular, we
let $\Sigma_l$ be the closure of $S_{k,l}$,
let $T_l$ be the (closed) tubular neighborhood of $\Sigma_l\sub R_k$ with 
the projection $p_l\mh T_l\to \Sigma_l$ and let $F_l$
be the subbundle of $W_k|_{T_l}$ with the isomorphism \eqref{1.10}.
Let $\zeta_l(t)$ be the standard extension of 
$(f_k^l)\sta(e_l(t))|_{\Sigma_l}$
to $T_l$. Note that $h_k|_{T_l}+\zeta_l(t)$ is a regular extension of 
$(f_k^l)\sta(\bh_l(t))|_{\Sigma_l}$. Because $\{\bh_l\}_{l<k}$
satisfies condition $\cH_k$, for $l<m<k$ and $x\in \Sigma_l\cap
\Sigma_{m}$ we have 
$(f_k^l)\sta(\bh_l(t))(x)=(f^{m}_k)\sta(\bh_{m}(t))(x)$.
Now let
\begin{equation*}
A_l=p_l\upmo \bl(f_k^l)\upmo(S\pri_l)\br-
\bigcup_{l<m<k} p_{m}\upmo\bl(f_k^{m})\upmo(\clo(R\pri_{m}))\br
\end{equation*}
and let
\begin{equation*}
B_l=\clo (R_{k,l}\pri)-
\bigcup_{k>m>l}(f_k^{m})\upmo(S_{m}\pri).
\end{equation*}
Note that $\{A_l\}_{l<k}$ covers $\Int\bl\cup_{l<k}
T_l\br$, that $B_l\subb A_l$ and that $\{B_l\}_{l<k}$ \
is a collection of compact subsets of $R_k$.
Now let $\eps>0$ be sufficiently small. We choose a collection of
non-negative smooth functions $\{\rho_l\}_{l<k}$ that obeys the 
requirement that $\supp(\rho_l)\subb\Int(A_l
\cap T_l\ueps)$, that $\rho_l\equiv1$ in a neighborhood of $B_l$
and that $\sum_{l<k}\rho_l\equiv1$ in a neighborhood of 
$\cup_{l<k}\clo(R_{k,l}\pri)$. This is possible because
the last set is compact and is contained in $\Int(
\cup_{l<k}A_l)$. We set
$$\zeta(t)=\sum_{l<k}\rho_l\cdot \zeta_l(t).
$$
Now we check that for each $l<k$ the section
$h_k+\zeta(t)$ is a regular extension of $(f_k^l)\sta(\bh_l(t))$ in
a neighborhood of $\clo(R_{k,l}\pri)$.
Let $x$ be any point in $\clo(R_{k,l}\pri)$. We first consider the case 
where $x$ is contained in $B_{m}$ for some $m\geq l$.
Let $y=f_k^{m}(x)$. Note that $y\in S_{m}\pri$.
Then restricting to a sufficiently small neighborhood
of $x$ the section $h_k+\zeta(t)$ is equal to $h_k+\zeta_{m}(t)$.
Since $h_k+\zeta(t)$ is a regular extension of
$(f_k^{m})\sta(\bh_{m}(t))$ near $x$ and
since $\bh_{m}(t)$ is a regular extension of $(f^l_{m})\sta
(\bh_l(t))$ in a neighborhood of $y\in S_{m}\pri$, $h_k+\zeta(t)$ is
a regular extension of $(f^l_k)\sta(\bh_l(t))$ near $x$. 
We next consider the case where $x$ is not contained in any of the 
$B_{m}$'s. Let $\Lambda$ be the set of all $m>l$ 
such that $x\in (f_k^{m})\upmo(S_{m}\pri)$. Then for any
$m<k$ that is not in $\Lambda$, $\rho_{m}\equiv 0$ 
in a neighborhood of $x$. Here we have used the fact
that $\Sigma_m\cap T_l\ueps\sub\Sigma_m\cap\Sigma_l$ for
$0<\eps\ll 1$.
On the other hand, by induction hypothesis for each 
$m\in\Lambda$ the section $h_k+\zeta_{m}(t)$ is a regular
extension of $(f_k^l)\sta(\bh_l(t))$ near $x$. Therefore
since $\sum_{m\in\Lambda}\rho_{m}\equiv1$ near $x$, in a small
neighborhood of $x$
$$h_k+\zeta(t)=\sum_{m\in\Lambda}\rho_{m}\cdot (h_k+\zeta_{m}(t))
$$
is also a regular extension of $(f_k^l)\sta(\bh_l(t))$.

Our last step is to extend $\zeta(t)$ to $R_k$. 
We let $e_k(t)$ be a smooth family of sufficiently small sections
of $W_k$ such that $e_k(0)\equiv0$, that the restriction
of $e_k(t)$ to a neighborhood of $\cup_{l<k}\clo\bl
(f_k^l)\upmo(R_l\pri)\br$ is $\zeta(t)$ and such that the section
$\bh_k(1)$ is transversal to the zero section in a
neighborhood of $\clo(R_k\pri)$ in $S_k$. The last condition is
possible because $h_k+\zeta(1)$ is transversal to the zero section
in a neighborhood of 
$\cup_{l<k}\clo(R_{k,l}\pri)$.
Therefore, by possibly shrinking $S_l\pri$ while still keeping
$R_l\pri\subb S_l\pri$ for $l<k$ if necessary, we can find an
$S_k\pri\subb S_k$ satisfying $R_k\pri\subb S_k\pri$
such that the induction hypothesis $\cH_k$ holds
for $\{e_l\}_{l\leq k}$ and $\{S_l\pri\}_{l\leq k}$, except possibly the 
third condition.

We now show that the third condition of $\cH_k$ holds as well.
We only need to check the inclusion \eqref{eq:1.12} for $l=k$.
First, by Lemma \ref{1.9} we can find an open $S\subb S_k$
such that $R_k\pri\subb S$ and that
\begin{equation}
\lab{eq:1.13}
\bh_k\upmo(0)\cap (\clo(S)-R\pri_k)\sub
\bl\bigcup_{i< k}(f^i_k)\upmo(R_i\pri)\br
\bigcup\bl\bigcup_{i> k}f^k_i(R_{i,k}\pri)\br.
\end{equation}
Now let 
$$D_1=\bh_k\upmo(0)\bigcap (\clo(S)-R_k\pri)\bigcap 
\bl\bigcup_{i< k}(f^i_k)\upmo(R_i\pri)\br
$$
and let
$$D_2=\bh_k\upmo(0)\bigcap \bl\clo(S)-R_k\pri\br\bigcap
\bl\bigcup_{i> k}f^k_i(R_{i,k} \pri)\br.
$$
Since $\bh_{k}(t)$ are small perturbations of $h_k$, 
we can assume that $\bh_k(t)$ are chosen so that
for any $t\in[0,1]$ the left hand side of \eqref{eq:1.13}
is contained in the union of neighborhood $V_1$ of $D_1$ and 
a neighborhood
$V_2$ of $D_2$. We remark that if we choose $\{e_l\}_{l\leq k}$ so that
their $C^2$-norms are sufficiently small, then we can make $V_1$ and $V_2$
arbitrary small.
Then by Lemma \ref{1.11} the vanishing locus of $\bh_k(t)$ 
inside $V_1$ is contained in $\cup_{i\leq k}(f^i_k)\upmo(S_i\pri)$. On the
other hand, since $\cup_{i\geq k}f_i^k(R_{i,k}\pri)$ is 
open, it contains $V_2$ since $D_2$ is compact and $V_2\supset D_2$
is sufficiently small. This proves the inclusion \eqref{eq:1.12}.

Therefore, by induction we have found $\{S_k\pri\}_{k\in\cL}$ and
$\{e_k(t)\}_{k\in\cL}$ that satisfy the condition $\cH_k$  
for $k=\#(\cL)+1$. Now let $\bg_l(t)=\bh_l(t)|_{R_l\pri}$.
Then $\g(t)=\{\bg_l(t)\}_{l\in\cL}$ satisfies the condition of the
proposition. Note that the left hand side of \eqref{eq:1.11}
is compact because it is contained in the union of compact
sets $\{\iota_k(\clo(R_k\pri))\}_{k\in\cL}$. This proves the
proposition.
\end{proof}

Let $\g(t)$ be the perturbation constructed by Proposition \ref{1.12}
with $\h=\bphi$. We define the Euler class of
$[\Phi\mh \bB\to\bE]$ to be the homology class in
$H\lsta(\bB;\QQ)$ represented by the current $[\g(1)\upmo(0)]$.
In the remainder of this section, we will sketch the argument that 
shows that this class is independent of the choice of the chart
$\R$ and the perturbation $\g$.

\begin{prop}
Let the notation be as before. Then the homology class $[\g(1)\upmo(0)]
\in H\lsta(\bB;\QQ)$ so constructed is independent of the choice of
perturbations.
\end{prop}

\begin{proof}
First, we show that if we choose two perturbations $\g_1(t)$ and
$\g_2(t)$ based on identical sub-atlas $\R\pri\subb\R$ as stated in 
Proposition \ref{1.12}, then we have $[\g_1(1)\upmo(0)]
=[\g_2(1)\upmo(0)]$. 
To prove this, all we need is to construct a family of perturbations
$\g_s(t)$, where $s\in [0,1]$, that satisfies conditions similar to that
of the perturbations constructed in Proposition \ref{1.12}. Since then we 
obtain a current
$$\bigcup_{s\in [0,1]}\iota(\g_s(1)\upmo(0))\times\{s\}\sub\bB\times [0,1]
$$
is a homotopy between the currents $\g_0(1)\upmo(0)_{\text{cur}}$
and $\g_1(1)\upmo(0)_{\text{cur}}$. The construction of
$\g_s(t)$ is parallel to the construction of $\g(t)$ in Proposition
\ref{1.12} by considering the data over $\{ R_k\times [0,1]\}_{k\in\cL}$.

Next, we show that the cycle $[\g(1)\upmo(0)]$ does not depend on the
choice of $\R_1\subb \R$. Let $\R_1\subb\R$ and $\R_2\subb\R$ be two
good sub-atlas and let $\g_1(t)$ and $\g_2(t)$ are two perturbations
subordinate to $\R_1$ and $\R_2$ respectively. Clearly, 
we can choose a sub-atlas $\R_0\subb\R$ such that $\R_1\sub\R_0$ and
$\R_2\sub\R_0$. Let $\g_0(t)$ be a perturbation given by Proposition
\ref{1.12} subordinate to $\R_0$. Then $\g_0(t)$ is also subordinate to
$\R_1$ and $\R_2$. Hence by the previous argument
$$[\g_1(1)\upmo(0)]=[\g_0(1)\upmo(0)]=[\g_2(1)\upmo(0)].
$$

It remains to show that the class $[\g(1)\upmo(0)]$ does not depend on the 
choice of the good atlas $\R$. For this, it suffices to show that 
for any two good atlas $\R$ and $\R\pri$ so that $\R$ is finer
than $\R\pri$, the respective perturbations $\g(t)$ and $\g\pri(t)$
gives rise to identical homology classes
$[\g(1)\upmo(0)]=[\g\pri(1)\upmo(0)]$. Let $\R=\{R_k\}_{k\in\cK}$ and
$\R\pri=\{ R_k\}_{k\in\cL}$, and let $U_k\sub\bB$ (resp. $U_l\sub\bB$)
be the open subsets so that $(R_k,W_k,\phi_k)$
(resp. $(R_l,W_l,\phi_l)$) are the smooth approximations of
$[\Phi\mh\bB\to\bE]$ over $U_k$ (resp. $U_l$) for $k\in\cK$
(resp. $l\in\cL$). As before, we denote 
$\iota_k\mh R_k\to U_k$ be the tautological map with $k\in\cK$ or
$k\in\cL$. Let $U_{kl}=U_k\cap U_l$.
We consider the good atlas $\R_0$ with charts
$R_{kl}=\iota_k\upmo(U_{kl})$, where $(k,l)\in\cK\times\cL$,
with bundles $W_{kl}=W_k|_{R_{kl}}$ and 
$\phi_{kl}$ the restriction of $\phi_k$,
where $R_{kl}$ is considered to be an open subset of $R_k$.
Using the extension technique in the proof of Proposition \ref{1.12},
we can construct a perturbation $\g_0(t)$ that is a regular
extension of $\g\pri(1)\upmo(0)$ under the obvious 
$\iota_l\upmo(U_{kl})\to R_{kl}$ and $W_l|_{\iota_l\upmo(U_{kl})}
\sub W_{kl}$. Therefore, $[\g_0(1)\upmo(0)]=[\g\pri(1)\upmo(0)]$.
On the other hand, since $W_{kl}=W_k|_{R_{kl}}$, $\g(t)$ induces a
perturbation $\g_0\pri(t)$ subordinate to $\R_0$. Hence,
$[\g_0(1)\upmo(0)]=[\g_0\pri(1)\upmo(0)]=[\g(1)\upmo(0)]$.
This proves the proposition.
\end{proof}

\section{Analytic charts}

The goal of this section is to construct a collection of local
smooth approximations $(\Lambda,V)$ so that the data
$\phi_V\mh R_V\to V|_{R_V}$ are analytic. Namely, $R_V$ are
complex manifolds, $V|_{R_V}$ are holomorphic vector bundles and 
$\phi_V$ are holomorphic
sections. In the next section
we will show that such $\phi_V$'s are Kuranishi
maps,  and hence the cones $\lim_{t\to\infty}\Gamma_{t\phi}$ are
the virtual cones constructed in \cite{LT1}.

We will use the standard notation in complex geometry in this section.
For instance, if $M$ is a complex manifold, we will denote by $T_xM$
the complex tangent space of $M$ at $x$ unless otherwise is mentioned. 
We will use complex dimension throughout this section, unless otherwise is
mentioned. Accordingly
the complex dimension of a set is half of its real dimension.
We will use the words analytic and holomorphic interchangably in this
section as well.

We begin with the construction of such local smooth approximations.
Let $w\in \bB$ be any point representing a holomorphic
stable map $f\mh\Sigma\to X$ with $n$-marked points. 
We pick a uniformizing chart 
$\Lambda=(\tilu,\bE\ltilu,\Phi\ltilu,G\ltilu)$
of $w$ over $U\sub\bB$ such that the elements of
$\tilu$ are stable maps $f_1\mh\Sigma_1\to X$ with (distinct)
$(n+k)$-marked points $\{x_i\}$ so that $\{f_1(x_m)\}_{m=n+1}^{n+k}$
are the $k$-distinct points of $f_1\upmo(H)$, where $H$ is a smooth complex
hypersurface of $X$ in general position of degree $k=[H]\cdot
[A]$ and $A=f_{\ast}([\Sigma])$,
and that the stable maps resulting from discarding the
last $k$ marked points of $f_1$ are in $U$. Here as
usual we assume that $U$ is sufficiently small
so that all stable maps in $U$ intersect $H$ transversally and 
positively. Note
that the later correspondence is the projection $\pi_U\mh
\util\to U$.
Let $\cY$ over $\tilde U$ be the universal (continuous) 
family of curves with $(n+k)$ marked
sections and let $\cF\mh\cY\to X$ be the
universal map. 

We let $\pi\mh \tilu\to\mgnk$ be the tautological map induced by
the family $\cY$ with its marked sections. Here
$\mgnk$ is the moduli space of $(n+k)$-pointed stable curves of
genus $g$. 
Without loss of generality, we can assume that no fibers of
$\cY$ with the marked points have non-trivial
automorphisms. It follows that $\mgnk$ is smooth near
$\pi(\tilu)$. As in section 1, we view $G\ltilu$ as a subgroup
of $S_{n+k}$. Then $G\ltilu$ acts on $\mgnk$ by permuting the
$(n+k)$-marked points of the curves in $\mgnk$, and the map
$\pi\mh\tilu\to\mgnk$ is $G\ltilu$-equivariant. Now let
$O\sub\mgnk$ be a smooth $G\ltilu$-invariant open neighborhood
of $\pi(\tilu)\sub\mgnk$ and let $p\mh\cX\to O$ be the universal family of 
stable curves over $O$ with $(n+k)$ marked sections
(In this section we will work with the analytic category unless otherwise
is mentioned). It follows that the $G\ltilu$-action on $O$
lifts to $\cX$ that permutes its marked sections.
For convenience, we
let $\cX\times_O\tilde U$ be the topological subspace of 
$\cX\times\tilde U$
that is the preimage of $\Gamma_{\pi}\sub O\times\tilde U$ under 
$\cX\times\tilde U\to O\times\tilde U$, where 
$\Gamma_{\pi}\sub O\times\tilde U$
is the graph of $\pi\mh\tilde U\to O$. Since no fibers of $\cY$
(with marked points) have non-trivial automorphisms, there is a 
canonical $G\ltilu$-equivariant isomorphism
$\cY\cong \cX\times_O\tilde U$ as family of pointed curves. 
Let $\pi_{\cX}$
and
$\pi\ltilu$ be the first and the second projection of $\cX\times_O\tilde U$.

Next, we let $(\cX_n,O_n;\Sigma,p_n,\varphi_n)$ be a semi-universal
family of the $n$-pointed curve $\Sigma$. Namely, $\cX_n$ is a 
(holomorphic) family of pointed prestable curves over 
the pointed smooth complex manifold $p_n\in O_n$ whose
dimension is equal to $\dim_{\CC}\Ext^1(\om\lsig(D),\cO\lsig)$,
where $D\sub\Sigma$ is the divisor of the $n$-marked points of 
$\Sigma$, $\varphi_n\mh \Sigma\to\cX_n|_{p_n}$ is an
isomorphism of $\Sigma$ with the fiber of $\cX_n$ over $p_n$ as
$n$-pointed curve, and the Kuranishi map
$T_{p_n}O_n\to \Ext^1(\om\lsig(D),\cO\lsig)$ of the family $\cX_n$
is an isomorphism. Note that $G\ltilu$ acts canonically on $\Sigma$.
For convenience, we let $\Pi_n(\cX)$ be the family of curves over $O$
that is derived from $\cX$ by discarding its last $k$-marked sections.
We now let $B=\pi_U\upmo(w)$ and
fix a $G\ltilu$-equivariant isomorphism
\begin{equation}
\coprod_{z\in B}\Pi_n(\cX)|_z\lra B\times \Sigma
\lab{eq:2.50}
\end{equation}
over $B$. 
Let $\aut_{p_n}(\cX_n)$ be
the group of those biholomorphisms of $\cX_n$ that keep the
fiber $\cX_n|_{p_n}$ invariant, that send fibers of
$\cX_n$ to fibers of $\cX_n$ and that fix the $n$-sections
of $\cX_n$. Possibly after shrinking $O_n$ if necessary,
we can assume that there is a homomorphism $\rho\mh G\ltilu\to
\aut_{p_n}(\cX_n)$ such that for any $\sigma\in G\ltilu$ the
$\rho(\sigma)$ action on $\cX_n|_{p_n}$ is exactly the $\sigma$
action on $\Sigma$ via the isomorphism $\varphi_n$. Finally, 
possibly after shrinking $U$ and $O$, we can
pick a $G\ltilu$-equivariant holomorphic
map $\varphi\mh O\to O_n$ such that $\varphi(B)=p_n$ and that
there is a $G\ltilu$-equivariant isomorphism of $n$-pointed curves
$\tilde{\varphi}\mh\cX\to O\times_{O_n}\cX_n$ that
extends the isomorphism \eqref{eq:2.50}.
We remark that the reason for doing this is to ensure that the
smooth approximation we are about to construct is $G\ltilu$-equivariant.

Next, we let $l$ be an integer to be specified later
and let $U_i\sub\Sigma$, $i=1,\ldots,l$,
be $l$ disjoint open disks away from the marked points and the nodal points
of $\Sigma$. We assume that $\cup_{i=1}^l U_i$ is
$G\ltilu$-invariant and that for any $\sigma\in G\ltilu$
whenever $\sigma(U_i)=U_i$ then $\sigma|_{U_i}={\mathbf 1}_{U_i}$.
By shrinking $U$, $O$ and $O_n$
if necessary, we can find disjoint open subsets
$\cU_{n,i}\sub\cX_n$ such that $\cup_{i=1}^l\cU_{n,i}$ is
$G\ltilu$-invariant, that $\cup_{i=1}^l\cU_{n,i}$ is $G\ltilu$-equivariantly
biholomorphic to $O\times\cup_{i=1}^l U_i$, that
$\cU_{n,i}\cap\Sigma=U_i$ and that the projections $\cU_{n,i}\to O$ induced
by the projection
$\cX\to O$ is the first projection of $O\times U_i$ ($=\cU_{n,i}$).
For convenience, for each $i$
we will fix a biholomorphism between $U_i$ and the
unit disk in $\CC$, and will denote by $U_i\uhf$ the open disk in 
$U_i$ of radius $1/2$. 
We let $\cU_i$ be the 
disjoint open subsets of $\cY$ defined by
$$\cU_i=\cU_{n,i}\times_{O_n}\tilu\sub\cX_n\times_{O_n}\tilu
\cong\cX\times_O\tilu\cong\cY.
$$
We will call $U_i$ and $\cU_i$ the distinguished open subsets
of $\Sigma$ and $\cY$ respectively.
Without loss of generality, we can assume that $\cup_{i=1}^l\cU_i$ is
disjoint from the $(n+k)$-sections of $\cY$. We also assume that there
are holomorphic coordinate charts $V_i\sub X$ so that 
$\cF(\cU_i)\sub V_i$.
We let $(w_{i,1},\ldots,w_{i,m})$, where $m=\dim
X$, be the coordinate variable of $V_i$
and let $\bv_i=\partial/\partial w_{i,1}$.
For each $i$ we pick a nontrivial $(0,1)$-form $\gamma_i$ on $U_i$ 
with $\supp(\gamma_i)\subb U_i\uhf$. We demand further that
if there is a $\sigma\in G\ltilu$ so that $\sigma(U_i)=U_j$ then 
$V_i=V_j$ as coordinate chart and $\sigma\sta(\gamma_i)=\gamma_j$.
We then let $\sigma_i$ be the $(0,1)$-form over $\cU_i$ with 
values in $\cF\sta(T_X)|_{\cU_i}$ that
is the product of the pull back of $\gamma_i$ via
$\cU_i\times_O\tilu\to U_i$ with $\cF\sta(\bv_i)|_{\cU_i}$,
and let $\tilde{\sigma}_i$ be the section over $\cY$ that is
the extension of ${\sigma}_i$ by zero. Obviously, $\tilde{\sigma}_i$ is 
a section of
$\bE\ltilu$, and $(\tilde{\sigma}_1,\ldots,\tilde{\sigma}_l)$ is
linearly independent fiberwise. Hence it spans a complex
subbundle of $\bE\ltilu$,
denoted by $V$. It follows from the construction that $V$ is 
$G\ltilu$-equivariant.

As in the previous section, we let $R=\Phi\ltilu\upmo(V)$,
let $W=V|_R$ and let $\phi\mh R\to W$ be the lifting of 
$\Phi\ltilu|_R\mh R\to \bE\ltilu|R$. The main task of this 
section is to show that
we can choose $U_i$, $\gamma_i$ and $V_i$ so that $R$ admits a 
canonical complex structure and that the section $\phi$ is holomorphic when 
$W$ is endowed with the holomorphic structure so that the
basis $\tilde{\sigma}_1|_R,\ldots,\tilde{\sigma}_l|_R$ is holomorphic.

To specify our choice of $U_i$, $\gamma_i$ and $V_i$, we need first to 
define the Dolbeault cohomology of holomorphic vector bundles over 
singular curves. Let $\cE$ be a locally free sheaf of $\cO\lsig$-modules
and let $E$ be the associated vector bundle, namely, $\cO\lsig(E)=\cE$.
We let $\om^0_{{\rm cpt}}(E)$ be the sheaf of smooth sections of $E$ that
are holomorphic in a neighborhood of $\sing(\Sigma)$ and let
$\om_{{\rm cpt}}^{0,1}(E)$ be the sheaf of smooth sections of $(0,1)$-forms
with values in $E$ that vanish in a neighborhood of $\sing(\Sigma)$.
Let
$$\dbar:
\Gamma(\om^0_{{\rm cpt}}(E))\lra \Gamma(
\om_{{\rm cpt}}^{0,1}(E))
$$
be the complex that send $\varphi\in \om^0_{{\rm cpt}}(E))$ to
$\dbar(\varphi)$. Since $\varphi$ is holomorphic near nodes of $\Sigma$,
$\dbar(\varphi)$ vanishes near nodes of $\Sigma$ as well, and hence
the above complex is well defined.
We define the Dolbeault cohomology
$H^0_{\dbar}(E)$ and $H^1_{\dbar}(E)$
to be the kernel and the cokernel of
$\dbar$.

\begin{lemm}
\lab{2.0}
Let $H^i(\cE)$ be the C\v{e}ch cohomology of the sheaf $\cE$.
Then there are canonical isomorphisms
$H^0_{\dbar}(E)\cong H^0(\cE)$ and $\Psi\mh
H^{0,1}_{\dbar}(E)\cong H^1(\cE)$.
\end{lemm}

\begin{proof}
The proof is identical to the proof of the classical result that the
Dolbeault cohomology is isomorphic to the C\v{e}ch cohomology 
for smooth complex manifolds. Obviously, $H^0_{\dbar}(E)$
is canonically isomorphic to $H^0(\cE)$. We now construct $\Psi$.
We first cover $\Sigma$ by open subsets $\{W_i\}$ so that the 
intersection of any of its subcollection is contractible. Now let
$\varphi$ be any global section in $\om^{0,1}_{{\rm cpt}}(E)$.
Then over each $W_i$ we can find a smooth function 
$\eta_i\in\Gamma_{W_i}(\om_{{\rm cpt}}^0(E))$ such that 
$\dbar\eta_i=\varphi|_{W_i}$. Clearly, the class in $H^1(\cE)$ 
represented by the cocycle $[\eta_{ij}]$, where
$\eta_{ij}=\eta_i|_{W_i\cap W_j}-
\eta_j|_{W_i\cap W_j}$, is independent of the choice of $\eta_i$,
and thus defines a homomorphism 
$\Gamma(\om_{{\rm cpt}}^{0,1}(E))\to
H^1(\cE)$. It is routine to check that it is surjctive and its 
kernel is exactly $\image(\overline{{\mathbf\partial}})$. 
Therefore, we have $H^{0,1}_{\dbar}(E)\cong H^1(\cE)$. Also, it
is direct to check that this isomorphism does not depend on the  choice
of the covering $\{W_i\}$. This proves the lemma.
\end{proof}

For any $z\in\util$, we denote by $\tilde{\sigma}_i(z)$ the restriction
of $\tilde{\sigma}_i$ to the fiber of $\cY$ over $z$.
We now choose the $l$ open disks $U_i\sub
\Sigma$, the $(0,1)$-forms $\gamma_i$ on 
$U_i$ and the coordinate charts $V_i\sub X$ such that
for any $\tilde w\in\pi_U\upmo(w)$ the collection
$\tilde{\sigma}_1(\tilde w),\ldots,
\tilde{\sigma}_l(\tilde w)$
spans $H^{0,1}_{\dbar}(f\sta T_X)$.
This is certainly
possible if we choose $l$ large
because the locus of $U_i$ are arbitrary as long as
they are away from the nodal points of $\Sigma$ and the marked
points, and the charts $V_i$ can also be chosen with a lot of
choice.
We fix once and for all such choices of $U_i$, $V_i$ and $\gamma_i$.
We then let $\cU_i\sub\cY$, $V\to\bE\ltilu$ and $R=\Phi\ltilu\upmo(V)$
be the objects constructed before according to this choice of
$U_i$, $\gamma_i$ and $V_i$. Let $\cY_R\to R$ be the restriction
to $R\sub\tilde U$ of the family $\cY\to\tilde U$
with the marked sections and let $F\mh \cY_R\to X$ be the 
associated map. 
We also fix a smooth function $\eta_i$ over $U_i$ so that 
$\dbar\eta_i=\gamma_i$. We next extend the collection $\{ U_i\}_{i=1}^l$
to an open covering $\{U_i\}_{i=1}^L$ so that
the intersection of any subcollection of $\{U_i\}$ are 
contractible, and that for any $i\leq l$ and $j\geq l+1$ the sets
$U_i\uhf$ and $U_j$ are disjoint. For convenience, we agree
that $\eta_j=0$ for $j>l$

 From now on, we will fix an $\tilde w\in R$ over $w$.

\begin{lemm}\lab{2.1}

There is a constant $A$ such that
for any C\v{e}ch 1-cocycle $[\tau_{ij}]$, where 
$\tau_{ij}\in\Gamma_{U_i\cap U_j}(f\sta \cT_X)$,
there are constants $a_i$ and holomorphic sections
$\zeta_i\in\Gamma_{U_i}(f\sta \cT_X)$ for $i=1,\ldots,L$ such that
$$(\zeta_j+a_j\eta_j)|_{U_j\cap U_i}-(\zeta_i+a_i\eta_i)|_{U_j\cap U_i}
=\tau_{ji}
$$
and 
$$\sum_{i=1}^L\bl\parallel\!\zeta_i\!\parallel_{L_2}+|a_i|\br\leq 
A\bl\sum_{i,j}^L\parallel\!\tau_{ij}\!\parallel_{L_2}\br.
$$
\end{lemm}

\begin{proof}
The existence of $\{a_i\}$ and $\{\zeta_i\}$ follows from
the fact that the images of $\tilde\sigma_1(\tilde w),
\ldots,\tilde\sigma_l(\tilde w)$ spans $H_{\dbar}^{0,1}(f\sta T_X)$
and that $H_{\dbar}^{0,1}(f\sta T_X)$ is isomorphic to $H^1(f\sta\cT_X)$.
The elliptic estimate is routine, using the harmonic theory on the
normalization of $\Sigma$. We will leave the details to the readers.
\end{proof}

We let $\Gamma(\om_{{\rm cpt}}^{0,1}(f\sta T_X))\shar$ be the quotient of
$\Gamma(\om_{{\rm cpt}}^{0,1}(f\sta T_X))$ by the linear span of 
$\tilde{\sigma}_1(\tilde w),\ldots,\tilde{\sigma}_l(\tilde w)$.
Because $\{\tilde\sigma_i(\tilde w)\}_{i=1}^L$
is invariant under the automorphism group of the stable
map $f$, $\Gamma(\om_{{\rm cpt}}^{0,1}(f\sta T_X))\shar$
is independent of the choice of $\tilde w\in\pi_U\upmo(w)$. We let
$${\dbar}\shar\mh \Gamma(\om_{{\rm cpt}}^0(f\sta T_X))\to
\Gamma(\om_{{\rm cpt}}^{0,1}(f\sta T_X))\shar
$$
be the induced complex. We define $H^0_{\dbar}(f\sta)\shar$
and $H^{0,1}_{\dbar}(f\sta T_X)\shar$ be the kernel and the
cokernel of the above complex. 

\begin{coro}\lab{2.2}
Let the notation be as before. Then
$H^{0,1}_{\dbar}(f\sta T_X)\shar=0$. Further,
the complex dimension of $H^0_{\dbar}(f\sta T_X)\shar$
is $\deg(f\sta T_X)+m(1-g)+l$.
\end{coro}

\begin{proof}
The vanishing of $H_{\dbar}^{0,1}(f\sta T_X)\shar$ follows from the 
surjectivity of ${\dbar}\shar$. The second part follows from
$$\begin{array}{ll}
&\dim H^0_{\dbar}(f\sta T_X)\shar-\dim H^{0,1}_{\dbar}(f\sta
T_X)\shar\\
=&\dim H^0_{\dbar}(f\sta T_X)-\dim H^{0,1}_{\dbar}(f\sta T_X)+l
=\chi(f\sta T_X)+l
\end{array}
$$
and the Riemann-Roch theorem.
\end{proof}

Next, we will describe the tangent space of $R$ at $\tilde w$.
By the smoothness result of \cite{LT2}, we know that $R$ is a 
smooth manifold of (complex) dimension $r\lexp$. 
As before, we let $D\sub \Sigma$ be the divisor of the first
$n$-marked points of $\tilde w$.
Since $f$ is holomorphic, $df\dual$ is a homomorphism
of sheaves $f\sta\om_X\to\om\lsig$.
We let 
$$\ddotw=[f\sta\Omega_X
\mapto{\alpha}\om_{\Sigma}(D)]
$$
be the induced complex indexed at $-1$ and $0$.
We will first define the extension space
$\extfc$ and then show that it is canonically isomorphic to 
$T_{\tilde w}R$. 
We begin with some more notations. Let $\cE$ be a
sheaf of $\cO\lsig$-modules that is locally free 
away from the nodal points of $\Sigma$. Then there is a
holomorphic vector bundle $E$ over $\Sigma^0$, where $\Sigma^0$
is the smooth locus of $\Sigma$, such that $\cO_{\Sigma^0}(E)
=\cE|_{\Sigma^0}$. 
We define $\cE\ua$ to be the sheaf so that the germs of
$\cE\ua$ at nodal points $p\in\Sigma$ (resp. smooth points
$p\in\Sigma^0$) are isomorphic to the
germs of $\cE$ at $p$ (resp. germs of $\om_{\Sigma^0}^0(E)$ at $p$).
The set $\extfc$ is the set of equivalence classes
of pairs $(v_1,v_2)$ as follows. The data $v_1$ is an element in
$\Ext^1(\om\lsig(D),\cO\lsig)$, which defines an exact sequence
\begin{equation}
\begin{CD}
0@>>>\cO\lsig@>{\varphi_1}>>
\cB@>{\varphi_2}>>\om\lsig(D)@>>> 0,
\end{CD}
\lab{eq:2.2}
\end{equation}
and equivalently a family of $n$-pointed nodal curves over 
$T=\spec
\CC[t]/(t^2)$, say $\cC_T$ with $n$-marked sections $\tilde{x}_i$
(See \cite[section 1]{LT1}).
Note that $\cB$ is locally free over $\Sigma^0$. 
The data $v_2$ is a homomorphism $f\sta\om_X\to \cB\ua$
such that, first of all, the diagram
\begin{equation}
\begin{CD}
@. @. f\sta\om_X @= f\sta\om_X\\
@. @. @Vv_2VV @V{df\dual}VV\\
0 @>>> \cO\lsig\ua @>{\varphi_1}>> \cB\ua @>{\varphi_2}>> 
\om\lsig(D)\ua @>>> 0
\end{CD}
\lab{eq:2.3}
\end{equation}
is  commutative, where the lower sequence
is induced by \eqref{eq:2.2}. Secondly, since $v_2$ is 
holomorphic near nodes of $\Sigma$,
the differential $\dbar v_2$ vanishes
near nodes of $\Sigma$, and since $df\dual$ is holomorphic, $\dbar v_2$
lifts to a global section $\beta$ of $\om\lcpt^{0,1}(f\sta T_X)$.
We require that there are constants $a_1,\ldots,a_l$ such that
$\beta=\sum a_i\tilde\sigma_i(\tilde w)$.

The equivalence relation of such pairs 
are the usual equivalence relation of the diagrams \eqref{eq:2.3}. 
Namely, two pairs
$(v_1,v_2)$ and $(v\pri_1,v_2\pri)$ with the associated data 
$\{\cB,\varphi_i\}$ and $\{\cB\pri,\varphi_i\pri\}$ are equivalent if
there is an isomorphism $\eta\mh\cB\to\cB\pri$ so that
$\eta\circ\varphi_1=\varphi_1\pri$, $\varphi_2=\varphi_2\pri\circ\eta$
and $\eta\circ v_2=v_2\pri$.

\begin{lemm}\lab{2.3}
Let the notation be as before. Then 
$\Ext^1\lsig(\ddotw,\cO\lsig\ua)\shar$ is canonically 
a complex vector space of
complex dimension $r\lexp$.
\end{lemm}

\begin{proof}
The fact that $\Ext^1\lsig(\ddotw,\cO\lsig\ua)\shar$ forms a complex 
vector space can be established using
the usual technique in homological algebra. For instance, if
$r\in \Ext^1\lsig(\ddotw,\cO\lsig\ua)\shar$ 
is represented by $\{\cB,\varphi_i,v_2\}$
shown in the diagram \eqref{eq:2.3}, then for any complex number 
$a$ the element $ar$ is
represented by the same diagram with $\varphi_1$ replaced by
$a\varphi_1$. We now prove that
\begin{equation}
\dim\Ext^1\lsig(\ddotw,\cO\lsig\ua)\shar=r\lexp.
\lab{eq:2.4}
\end{equation}
Clearly, the following familiar sequence is still exact
in this case:
$$
\Ext^0(\cD^{\bullet}_f,\cO\lsig\ua)\lra \Ext^0(\om\lsig(D),\cO\lsig)
\lra H^0_{\dbar}(f\sta T_X)\shar\qquad\qquad\qquad
$$
$$
\qquad
\lra \Ext^1(\om\lsig(D),\cO\lsig)\lra \Ext^1(\om\lsig(D),\cO\lsig)\shar
\lra H^{0,1}_{\dbar}(f\sta T_X)\shar.
$$
Since $f$ is stable, $\Ext^0(\cD^{\bullet}_f,\cO\lsig\ua)=0$. 
Hence \eqref{eq:2.4} follows from Corollary \ref{2.2} and the Riemann-Roch
theorem. This proves the lemma.
\end{proof}

Recall that
should $R$ be a scheme then the Zariski tangent space of $T_{\tilde w}R$
would be the space of morphisms $\spec \CC[t]/(t^2)\to R$
that send their only closed points to $\tilde w$ modulo certain equivalence
relation. In the following, 
we will imitate this construction and construct the space
of pre-$\CC$-tangents of $R$ at $\tilde w$.
We still denote by $U_1,\ldots, U_l$ the $l$-distinguished open subsets of
$\Sigma$ and let
$\{U_i\}_{i=1}^L$ be an extension of $\{U_i\}_{i=1}^l$
to an open covering of $\Sigma$ such that the intersection of any 
of its subcollection are contractible. 
Without loss of generality, we assume $U_j\cap U_i\uhf=\emptyset$ for
$j>l$ and $i\leq l$. We also assume that there are coordinate charts $V_i$ of
$X$ such that $f(U_i)\sub V_i$. 
By abuse of notation, we will fix the embedding $V_i\sub \CC^m$
and view any map to $ V_i$ as a map to $\CC^m$.
We let $\iota_i\mh V_i\to X$ be the 
tautological inclusion and let 
$$g_{ij}\mh \iota_j\upmo(\iota_i(V_i))\to \iota_i\upmo(\iota_j(V_j))
\sub V_i\sub \CC^m
$$
be the transition functions of $X$.

We define a pre-$\CC$-tangent $\xi$ of $R$ at $\tilde w$ 
to be a collection of data as follows: First, there is a flat
analytic family of $n$-pointed pre-stable curves $C_T$ over an open
neighborhood $T$ of $0\in \CC$ such that
the fiber of $C_T$ over $0$, denoted by $C_0$, is isomorphic to $\Sigma$ 
as $n$-pointed curve; Secondly, there is
an open covering  $\{\tilde U_i\}_{i=1}^L$ of
$C_T$ such that $\tilde U_i\cap C_0=U_i$, and that
for each $i\leq l$, there is a biholomorphism $\tilde U_i\cong U_i\times T$ 
such that its restriction to $U_i=\tilde U_i\cap
C_0$ is compatible to the identity map of $U_i$;
Thirdly, there is a collection of smooth maps 
$\tilde f_i\mh \tilde U_i\to V_i$ such that for $i>l$, 
all $\tilde f_i$ are
holomorphic and that for each
$i\leq l$ we have 
${\dbar}_0(\tilde f_i)=0$ and 
\begin{equation}
{\dbar}_i(\tilde f_i)=
\pi_T\sta\varphi_i \cdot\pi_{U_i}\sta(\gamma_i\cdot f\sta(\bv_i)),
\lab{eq:2.22}
\end{equation}
where $\pi_{U_i}$ and $\pi_T$ are the first and the second projection of
$U_i\times T$, $\varphi_i$ are holomorphic functions over $T$
and $\dbar_0$ (resp. $\dbar_i$) is the $\dbar$-differential with respect to
the holomorphic variable of $T$ (resp. $U_i$) using
$\tilde U_i\cong U_i\times T$
and the $\gamma_i$ and $\bv_i$ are
the $(0,1)$-form and the vector field chosen before;
Forthly, if we let $z_0$ be the holomorphic variable of $\CC\supset T$, then
we require that
\begin{equation}
\lab{eq:2.27}
\tilde f_{ji}=\tilde f_i- g_{ij}\circ \tilde f_j:
\tilde U_{ij}\lra \CC^m,
\end{equation}
where $\tilde U_{ij}$ is a neighborhood of $U_i\cap U_j$ in
$\tilde U_i\cap \tilde U_j$ over which $\tilde f_{ji}$
is well-defined, is divisible by $z^2_0$ 
(Namely, $\tilde f_{ji}$ has the form
$\pi_T\sta(z^2_0)\cdot h_{ji}$ for some smooth function 
$h_{ji}\mh \tilde U_{ij}\to \CC^m$).
Intuitively, a pre-$\CC$-tangent is 
a scheme analogue of a morphism $\spec \CC[t]/(t^2)\to R$
should $R$ be a scheme. We denote the set of all pre-$\CC$-tangents by
$T_{\tilde w}\upre R$. Note that $T_{\tilde w}\upre R$ is merely a
collection of all pre-$\CC$-tangents.

We next define a canonical map
\begin{equation}
T_{\tilde w}\upre R\lra \Ext^1(\ddotw,\cO\lsig)\shar.
\lab{eq:2.26}
\end{equation}
Let $\xi$ be any pre-$\CC$-tangent given by the data above. By
the theory of deformation of $n$-pointed curves
\cite[section 1]{LT1}, the analytic family $C_T$ 
defines canonically an exact sequence
\begin{equation}
0\lra\cO\lsig\lra \cB\lra \om\lsig(D)\lra 0,
\lab{eq:2.23}
\end{equation}
associated to an extension class $v_1(\xi)\in\Ext^1(\om\lsig(D),\cO\lsig)$,
where away from the nodes of $\Sigma$ and the suport of $D$ the sheaf
$\cB$ is canonically isomorphic to $\om_{C_T}\otimes_{\cO_{C_T}}\cO_{C_0}$.
Because $\tilde f_i\mh \tilde U_i\to V_i$ are
holomorphic for $i>l$, it follows from \cite{LT1} that there is
a canonical homomorphism of sheaves $u_i\mh f\sta\om_X|_{U_i}
\to \cB|_{U_i}$ such that
\begin{equation}
\begin{CD}
@. @. f\sta\om_X|_{U_i} @= f\sta\om_X|_{U_i}\\
@. @. @VV{u_i}V @VV{df\dual|_{U_i}}V\\
0 @>>> \cO_{U_i}\ua @>>> \cB\ua|_{U_i} @>>> \om_{U_i}(D)\ua
@>>> 0\\
\end{CD}
\end{equation}
is commutative, where the lower sequence is induced by \eqref{eq:2.23}.
Indeed, at smooth point $p\in U_i$ away from the support of $D$ the dual of
$u_i\otimes k(p)$ is the differential
$$d\tilde f_i(p) : T_p C_T=\cB\dual\otimes k(p)\lra f\sta T_X|_p.
$$
Note that by our choice of $U_i$, for $i\leq l$ the distinguished
open subsets $U_i$ are disjoint from the support of the 
$(n+k)$-marked points of $\tilde w$. Hence $\cB\ua|_{U_i}$ are canonically
isomorphic to $\om_{U_i}^0(TC_T|_{U_i})$, and
the dual of $d\tilde f_i$ define  canonical homomorphisms $u_i\mh 
f\sta\om_X|_{U_i}\to\cB|_{U_i}$ that make the above diagrams commutative.
Because of the condition \eqref{eq:2.22}, the lift
of $\dbar u_i$ is a constant multiple of $\tilde\sigma_i(\tilde w)|_{U_i}$.
Further, because of the condition that $\tilde f_{ji}$ is divisible by $z_0^2$,
the collection $\{ u_i\}_{i=1}^L$
patch together to form a homomorphism $v_2(\xi)\mh f\sta\om_X\to\cB\ua$
that makes the diagram \eqref{eq:2.3} commutative. Hence $(v_1(\xi),v_2(\xi))$
defines an element in $\extfc$, which is defined to be the image of $\xi$.

We remark that in this construction we have only used the fact
that the stable map associated to $\tilde w$ is holomorphic, that
the domain $\Sigma$ of $\tilde w$ has $l$ distinguished open
subsets $U_i$ with $(0,1)$-forms $\tilde\sigma_i(\tilde w)$.
Because for any $z\in R$ its domain $\Sigma_z$ also has $l$ distinguished
open subsets, namely $\cU_i\cap \Sigma_z\cong U_i$, and the forms
$\tilde\sigma(z)$, we can define the extension group
$\Ext^1(\dbul_z,\cO_{\Sigma_z})\shar$, the space of pre-$\CC$-tangents
of $R$ at $z$ and the analogut canonical map as in \eqref{eq:2.26}
if the map $f_z$ of $z$ is holomorphic.

To justify our choice of $\extfc$, we will construct, to each $v\in\extfc$,
a pre-$\CC$-tanget $\xi^v\in T_{\tilde w}\upre R$ so that the
image of $\xi^v$ under \eqref{eq:2.26} is $v$.
Let $v=(v_1,v_2)$ be any element in $\extfc$ defined by the diagram
\eqref{eq:2.3}. Let $T\sub \CC$ be a neighborhood of $0$ and 
let $C_T$ be an analytic family of $n$-pointed curves so that $C_0\cong \Sigma$
and the Kuranishi map $T_0\CC\to\Ext^1(\om\lsig(D),\cO\lsig)$ send
$1$ to $v_1$. For instance, we can take $C_T$ be the
pull back of $\cX_n$ via an analytic map $(T,0)\to (O_n,p_n)$. We let
$\{U_i\}_{i=1}^L$ be a covering of
$\Sigma$ as before and let $\{\tilde U_i\}_{i=1}^L$ be a covering of $C_T$
that are the pull back of $\cU_{n,i}$. Note that for $i\leq l$, they 
come with biholomorphisms $\tilde U_i\cong U_i\times T$. Let
$V_i$ be open charts of $X$ as before with $f(U_i)\subb V_i$. For $i>l$,
since the restriction of \eqref{eq:2.3} to $U_i$ is analytic, we can
find analytic $\tilde f_i\mh \tilde U_i\to V_i$, possibly after 
shrinking $T$ if necessary, such that $\tilde f_i$ are related to 
$v_2|_{U_i}$ as to how $u_i$ are related to $v_2(\xi)|_{U_i}$ before.
By analytic analogue of deformation theory (see \cite{LT1})
such $\tilde f_i$ do exist. For $i\leq l$, since $U_i$ are smooth
and $\cB\ua|_{U_i}$ are the sheaves $\om_{U_i}^0(T\sta C_T|_{U_i})$,
we simply let $\tilde f_i$ be smooth so that in addition to $\tilde
f_i$ satisfying the condition on pre-$\CC$-tangents we require that
$v_2|_{U_i}$ coincide with the dual of $d\tilde f_1|_{U_i}$. Note that
$(C_T,\{\tilde f_i\})$ will be a pre-$\CC$-tangent if $\tilde f_{ji}$
in \eqref{eq:2.27}  is divisible by $\pi_T\sta(z_0^2)$.
But this is true because for any $p\in U_i\cap U_j$, the differential
$d\tilde f_i(p)$ and $d\tilde f_j(p)$ from $T_p C_T$ to $T_{f(p)}X$ 
are identical. We let the so constructed pre-$\CC$-tangent be $\xi^v$.
Of course $\xi^v$ are not unique. It is obvious from the construction that
the image of $\xi^v$ under \eqref{eq:2.26} is $v$.
We remark that it follows from the construction that for any complex
number $c\ne 0$ the pull back of $(C_T,\{\tilde f_i\})$ under $L_c
\mh \CC\to \CC$ defined by $L_c(z_0)=cz_0$ is a pre-$\CC$-tangent,
say $\xi^{cv}$, whose image under \eqref{eq:2.26} is $cv$.

We next construct a holomorphic coordinate chart of $R$ at $\tilde w$. 
Let $r=\dim R$, which is $r\lexp+l=\dim\extfc$.
We fix a $\CC$-isomorphism $T_0\CC^r\cong\extfc$. Composed with
the canonical
$$\extfc\to \Ext^1(\om\lsig(D),\cO\lsig),
$$
we obtain
\begin{equation}
T_0\CC^r\lra \Ext^1(\om\lsig(D),\cO\lsig).
\lab{eq:2.20}
\end{equation}
Let $\cX_n$ over $O_n$ be the semi-universal family of
the $n$-pointed curve $\Sigma$ given before. We let $S$ be a
neighborhood of $0\in\CC^r$ and let $\varphi\mh S\to O_n$
be a holomorphic map with $\varphi(0)=0$ such that
$$d\varphi(0): T_0S\equiv T_0\CC^r\lra T_{p_n}O_n\cong
\Ext^1(\om\lsig(D),\cO\lsig)
$$
is the homomorphism \eqref{eq:2.20}. We let $\pi_S\mh C_S\to S$ be the family
of $n$-pointed curves over $S$ that is the pull back of $\cX_n$. Note that
$C_S|_0$, denoted by $C_0$, is canonically isomorphic to $\Sigma$.

We keep the open covering $\{U_i\}^L_{i=1}$ of $\Sigma$ ($\cong C_0$)
chosen before. We let $\{W_i\}_{i=1}^L$ be an open covering
of a neighborhood of $C_0\sub C_S$ so that $W_i\cap C_0=U_i$.
For $i\leq l$, we let $W_i$ be the pull back of $\cU_{n,i}\sub \cX_n$. For
$i>l$ and $U_i$ smooth, we choose $W_i$ so that there is a 
holomorphic map $\pi_i\mh W_i\to U_i$ so that the restriction of
$\pi_i$ to $U_i$ is the identity map. For $i>l$ and $U_i$ contains
a nodal point, we assume that $W_i$ is biholomorphic to
the unit ball in $\CC^{r+1}$ so that $U_i\sub W_i$ is defined by
$w_1w_2=0$ and $w_i=0$ for $i\geq 3$,
where $(w_i)$ are the coordinate variables of $\CC^{r+1}$,
and the restriction
of $\pi_S$ to $W_i$ is given by
$$z_1=w_1w_2, \ z_2=w_3,\ldots,z_r=w_{r+1},
$$
where $(z_i)$ are the coordinate variables of $\CC^r$. The upshot
of this is that if $h$ is a holomorphic function on $U_i$,
then we can extend it canonically to $W_i$ as follows. In case 
$U_i$ is smooth, then the extension of $h$ is the 
composite of $W_i\to U_i$ with $h$;  In case
$U_i$ is singular, then $\varphi$ has a unique expression 
$$a+w_1h_1(w_1)+w_2h_2(w_2),
$$
where $a\in\CC$ and $h_1, h_2$ are holomophic.
We then let its extension be the holomorphic function on
$W_i$ that has the same expression.

We fix the choice of $\{U_i\}$ and $\{W_i\}$. Without loss of generality,
we can assume that there are coordinate charts $V_i\sub X$
so that $f(U_i)\subb V_i$. Of course, for $i\leq l$ the charts $V_i$ are
the charts we have chosen before. 
Our construction of the local holomorphic chart of $R$ is
parallel to the original construction of Kodaira-Spencer of semi-universal
family of deformation of holomorphic structures without obstructions.
To begin with, possibly after shrinking $W_i$ if necessary we
can assume that the maps $f|_{U_i}\mh U_i\to V_i$ can be extended
to a holomorphic $F_{0,i}\mh W_i\to V_i$ (Recall $f$ is holomorphic).
We now let $\cA(W_i,V_i)$ be the space of smooth maps 
from $W_i$ to $\CC^m$ defined as follows. 
If $i>l$, then $\cA(W_i,V_i)$ consists of holomorphic maps
from $W_i$ to $\CC^m$; If $i\leq l$, then using the isomorphism
$W_i\cong U_i\times S$ and holomorphic coordinate $z=(z_i)$ of
$S$ and holomorphic coordinate $\xi$ of $U_i$, any smooth function
$\varphi\mh W_i\to \CC^m$ can be expressed in terms of its
$m$ components $\varphi_j(z, \xi)$, $j=1,\ldots,m$. We 
define $\cA(W_i,V_i)$ to be the set of
those smooth maps $\varphi\mh
W_i\to\CC^m$ so that 
$$\left\{
\begin{array}{l}
\dbar_{z_k}\varphi_j=0\quad \text{for}\quad k=1,\ldots,r \and j=1,\ldots,m;\\
\dbar_{\xi}\varphi_j=0\quad \text{for}\quad j\geq 2
\and \dbar_{\xi}\varphi_1=c\sigma_i\pri\quad \text{for some}\quad
c\in\CC,
\end{array}
\right.
$$
where $\sigma\pri_i$ is a $(0,1)$-form taking values in $\varphi\sta
\CC^n$ corresponding to the form $\sigma_i$ using the canonical
embedding $V_i\sub\CC^n$.
Note that $\cA(W_i,V_i)$ are $\cO_S$-modules. In
particular, if we let $\cI\sub \cO_S $ be the ideal sheaf of $0\in S$,
then we denote by $\cI^q \cA(W_i,V_i)$ the image of
$\cI^q\otimes_{\cO_S }\cA(W_i,V_i)$ in $\cA(W_i,V_i)$.

In the following, we will construct
a sequence of maps $F_{s,i}\in\cA(W_i, V_i)$ indexed by
$s\geq 1$ and $1\leq i\leq L$ of which the following holds:

\noindent
1. For each $i$, $F_{s+1,i}-F_{s,i}\in \cI^{s}\cA(W_i, V_i)$;

\noindent
2. The restrictions $F_{1,i}|_{U_i}\mh U_i\to\CC^m$ factor through 
$V_i\sub\CC^m$ and $\iota_i\circ (F_{1,i}|_{U_i})\mh U_i\to X$ 
is identical to $f|_{U_i}\mh U_i\to X$;

\noindent
3. In a neighborhood $W_{ij}$ of $U_i\cap U_j$ in $W_i\cap W_j$
over which the map
\begin{equation}
 F_{s,ij}=g_{ij}\circ F_{s,j}- F_{s,i}\mh W_{ij}\to \CC^m
\lab{eq:2.21}
\end{equation}
is well defined,
$F_{s,ij}\in \cI^{s}\cH(W_{ij},\CC^m)$, where $\cH(W_{ij},\CC^m)$
is the $\cO_S $-module of holomorphic maps from $W_{ij}$ to $\CC^m$;

\noindent
4. For any vector $\eta\in\CC^r$, we let $L_{\eta}\mh\CC\to\CC^r$ be the
unique $\CC$-linear map so that $L_{\eta}(1)=\eta$, and let $\eta\upre$ be 
the pre-$\CC$-tangent associated to the pull back of $(C_S,\{ F_{2,i}\})$
under $L_{\eta}$. Using the standard isomorphism $T_0 S\equiv T_0\CC^r\cong
\CC^r$, we obtain a map
\begin{equation}
T_0S\lra \extfc
\end{equation}
that send $\eta\in T_0 S$ to the image of $\eta\upre$ under \eqref{eq:2.26}.
We require that this map is the isomorphism \eqref{eq:2.20}.

For $s=1$ we simply let $F_{1,i}$ be the standard extension
of $f|_{U_i}\mh U_i\to V_i$ to $W_i\to \CC^m$.
We now show that we can construct $\{ F_{2,i}\}$ as required.
We let $\pi_1$ and $\pi_2$ be the first and the
second projection of $\CC^r\times \Sigma$, where
we view $\CC^r$ as the total space of $\extfc$. It follows
from the definition of the extension group that there is a
universal diagram
\begin{equation}
\begin{CD}
@.@.\pi_2\sta f\sta\om_X @= \pi_2\sta f\sta\om_X\\
@.@. @VV{\cV_2}V @VV{\pi_2\sta(df\dual)}V \\
0 @>>> \pi_2\sta\cO\lsig\ua @>>> {\tilde\cB}\ua 
@>>> \pi_2\sta\om\lsig(D)\ua @>>> 0
\end{CD}
\lab{eq:2.28}
\end{equation}
such that its restriction to fibers of $\CC^r\times\Sigma$ over
$\xi\in\CC^r$ are the diagrams \eqref{eq:2.3} associated to $\xi
\in\extfc$. By deformation theory of pointed curves, for any smooth
point $p\in\Sigma$ the vector space $\tilde\cB\otimes k(p)$ is
canonically isomorphic to the cotangent space $T_p\sta C_S$.
By applying the construction of $\xi^v\in T_{\tilde w}\upre
R$ from $v\in\extfc$ to the family version, we can construct the
family $\{ F_{2,i}\}$ as required. We will leave the details to the readers.

Now we show that we can successively construct $F_{s,i}$ 
that satisfies the four conditions above.
Assume that for some $s\geq 2$ we have constructed 
$\{F_{s,i}\}$ that satisfies the four conditions above. 
Let $W_{ij}$ be the neighborhood of $U_{ij}=U_i\cap U_j\sub
W_i\cap W_j$ so that \eqref{eq:2.21} is well-defined. Then by the
condition 3 above, $F_{s,ij}\in \cI^{s} \cH(W_{ij},\CC^m)$.
Let $I=(i_1,\ldots,i_r)$ be any length $s$ mulptiple index, 
namely, $i_j\geq 0$ and $\sum i_j=s$. 
As usual, we will denote by $\partial^I$ the
symbol $\partial^{i_1}/\partial z_1^{i_1}\cdots\partial^{i_r}/\partial
z_r^{i_r}$
and by $z^I$ the term $z_1^{i_1}\cdots z_r^{i_r}$.
Then because of the
condition 3 above,
$\varphi_{I,ij}={\partial^I} F_{s,ij}|_{U_{ij}}$
is a holomorphic section of $f\sta T_X|_{U_{ij}}$
using the standard isomorphism 
$$TX|_{V_i}\cong
TV_i\cong V_i\times\CC^m,
$$
and the collection
$[\varphi_{I,ij}]$ defines a C\v{e}ch 1-cocycle of 
$f\sta\cT_X$. We let $\{\phi_{I,i}\}$, where $\phi_{I,i}=\zeta_i+a_i\eta_i$,
be the collection provided by Lemma \ref{2.1}. 
Using the standard isomorphism $TX|_{V_i}\cong V_i\times\CC^m$,
we can view $\phi_{I,i}$ as a map $V_i\to \CC^m$. We let $\tilde\phi_{I,i}\mh
W_i\to\CC^m$ be the standard extension of $\phi_{I,i}$ and let $G_{I,i}=
\pi_S\sta(z^I)\tilde\phi_{I,i}$. Clearly, $\partial^I G_{I,i}=\phi_{I,i}$.
Now we let 
$$F_{s+1,i}=F_{s,i}+\sum_{\ell(I)=s}  G_{I,i}.
$$
It is direct to check that the collection $\{F_{s+1,i}\}$ 
satisfies the condition 1-4 before.
Finally, by the estimate in Lemma \ref{2.1}, there is a 
neighborhood of $U_i\sub W_i$, say $W_i^0$, such that $\lim_s F_{s,i}$
converges over $W_i^0$. Let
$F_{\infty,i}$ be its limit.
Because $f(U_i)\subb V_i$, there is a neighborhood 
$\tilde{W}_i$ of $U_i\sub W_i^0$ such that $F_{\infty,i}(\tilde{W}_i)
\sub V_i\sub\CC^m$. It follows that we can find a neighborhood
$S^0\sub S$ of $0\in S$ such that $\pi_S\upmo(S^0)\sub \cup \tilde{W}_i$.
Finally, because $F_{\infty,i}$ is analytic near $U_i$ for $i>l$
and is analytic in $S$ direction using $W_i\cong U_i\times S$ otherwise,
the condition 3 implies that
the collection $F_{\infty,i}|_{{W}_i\cap\pi_S\upmo(S^0)}$
defines a map 
$$F_S: \pi_S\upmo(S^0)\lra X.
$$ 
Clearly, $F_S$ is holomorphic away from the union of $W_1,\ldots, W_l$.
Further, for each $i\leq l$ if we let $\xi_i$ be a holomorphic
variable of $U_i$ and let $\pi_{U_i}$ and $\pi_{S^0}$ be the 
first and the second projection of 
$W_i\cap \pi_S\upmo(S^0)\cong U_i\times S^0$, then
\begin{equation}
\frac{\partial}{\partial\bar\xi_i} F_S|_{W_i\cap\pi_S\upmo(S^0)}d\bar\xi_i=
\pi_{S^0}\sta(\varphi_i)\pi_{U_i}\sta(\gamma_i) F_S\sta(\bv_i)|_{W_i}
\lab{eq:2.30}
\end{equation}
where $\varphi_i$ is a holomorphic
function over $S^0$.
Finally, we let $Z$ be the subset of 
$$\pi_S\upmo(S^0)\times_S\cdots\times_S\pi\upmo_S(S^0)\qquad
(\text{$k$ times})
$$
consisting of $(s;x_{n+1},\ldots,x_{n+k})$ such that $s\in S^0$
and that  
$ x_{n+1},\ldots,x_{n+k}$ are distinct points in $\pi_S\upmo(s)$ 
that lie in $F_s\upmo(H)$. Note that if we choose $U$ to be small
enough, then $F_s\upmo(H)$ has exactly $k$ points. 
Let $C_Z$ be the family of $(n+k)$-pointed curves over $Z$ so that
its domain is the pull back of $C_S$ via $Z\to S$, its first
$n$-marked sections is the pull back of the $n$-marked sections of $C_S$
and its last $k$-sections of the fiber of $C_Z$ over $(s;x_{n+1},\ldots,
x_{n+k})$ is $x_{n+1},\ldots,x_{n+k}$.
Coupled with the pull back of $F_S$, say $F_Z\mh C_Z\to X$,
we obtain a family of stable (continuous) maps from
$(n+k)$-pointed curves to $X$. Let $\eta\mh Z\to  U$ 
be the tautological map.

We claim that $\eta(Z)\sub R$. Indeed, let $z\in Z$ be any
point and let $C_z$ be the domain of $z$. It follows from
our construction that $C_z$ has $l$ distinguished open subsets,
denoted by $U_1,\ldots, U_l$, such that $f_z=F_Z|_{C_z}$
is holomorphic away from $\cup_{i=1}^l U_i$ and
$\dbar f_z|_{U_i}$ is a constant multiple of $\gamma_i\cdot
f_z\sta(\bv_i)$. Hence the value of the section $\Phi_{\tilde U}\mh
\tilde U\to \bE_{\tilde U}$ at $\eta(z)$ is contained in the
subspace $V|_{\eta(z)}\sub \bE_{\tilde U}|_{\eta(z)}$.
This shows that $\eta(z)\in R$.

\begin{prop}
The induced map $\eta\mh Z\to R$ is a local diffeomorphism near
those $z\in R$ whose associated map $f_z\mh C_z\to X$ are
holomorphic.
\end{prop}

\begin{proof}
This follows immediately from the proof of the basic Lemma in
\cite{LT2}. We will omit the details here.
\end{proof}

By shrinking $S^0$ if necesary, we can assume that
$\eta\mh Z\to R$ is a local diffeomorphism. We can further assume that
there is an open subset $Z\pri\sub Z$ containing $\tilde w$
such that $\eta\pri=\eta|_{Z\pri}\mh Z\pri\to R$ is one-to-one
and the image $\eta(Z\pri)\sub R$ is invariant under $G\ltilu$.
$\eta\pri\mh Z\pri\to R$ is the analytic coordinate of $\tilde w\in R$
we want. For convenience,
we will view $Z\pri$ as an open subset of $R$.

\begin{prop}
Let $V\pri$ be the restriction of $W$ to $Z\pri$ endowed with the
holomorphic structure so that $\tilde\sigma_1|_{Z\pri},\ldots,
\tilde\sigma_l|_{Z\pri}$ is a holomorphic frame. Then $\phi\pri
\equiv \phi_V|_{Z\pri}
\mh Z\pri\to V\pri$ is holomorphic.
\end{prop}

\begin{proof}
This follows immediately from \eqref{eq:2.30}.
\end{proof}

Let ${\phi\pri}\upmo(0)$ be any point and let $f_z\mh C_z\to X$ be the
associated (analytic) stable map with $D_z$ the divisor of its
first $n$-marked points. Then there is a caonical exact sequence of 
vector spaces
$$
\Ext^1(\Om_{C_z}(D_z),\cO_{C_z})\lra
H^1(f_z\sta\cT_X)\lra \Ext^2(\dbul_z,\cO_{C_z})\lra 0
$$
induced by the short exact sequence of complexes
$$
0\lra [0\to\Om_{C_z}(D_z)]\lra
[f_z\sta\Om_{X}\to\Om_{C_z}(D_z)]
\lra[f_z\sta\Om_X\to 0]\lra 0.
$$
Similarly, the differentil $d\phi_V(z)\mh T_z R\to W_z$ induces
an exact sequence of vector spaces
$$
\Ext^1(\dbul_z,\cO_{C_z}\ua)\sha\mapright{d\phi_V(z)} W|_z
\lra \coker(d\phi_V(z))\lra 0.
$$
Note that there are canonical homomorphisms $\Ext^1(\dbul_z,\cO_{C_z}
\ua)\sha\to\Ext^1(\dbul_z,\cO_{C_z})$ and $W|_z\to H^{0,1}_{\dbar}(
f_z\sta T_X)\cong H^1(f_z\sta\cT_X)$.

\begin{lemm}
\lab{2.22}
There is a canonical isomorphism $\xi$ (as shown below)
that fits into the diagram
\begin{equation*}
\begin{CD}
\Ext^1(\dbul_z,\cO_{C_z})\sha @>{d\phi_V(z)}>> W|_z
@>>> \coker(d\phi_V(z)) @>>> 0\\
@VVV @VVV @VV{\xi}V\\
\Ext^1(\Om_{C_z}(D_z),\cO_{C_z}) @>>> H^1(f_z\sta\cT_X) @>>>
\Ext^2(\dbul_z,\cO_{C_z}) @>>> 0.
\end{CD}
\end{equation*}
\end{lemm}

\begin{proof}
This is obvious and will be left to the readers.
\end{proof}

\section{The proof of the comparison theorem}

In this section, we will prove that the algebraic and the symplectic
construction of GW-invariants yield identical invariants.

We will work with the category of algebraic schemes 
as well as the category of 
analytic schemes. Specifically, we will use the words schemes, 
morphisms and \`etale neighborhoods  
to mean the corresponding objects in algebraic category and use the
word analytic maps and open subsets to mean the corresponding objects in
analytic category. As before, the words analytic and holomorphic are
interchangable. Also, we will use $\cO_S$ to mean the
sheaf of algebraic sections or the sheaf of analytic sections depending on
whether $S$ is an algebraic scheme or an analytic scheme. 
We will continue to use the complex dimension through out this section.

We now clarify our usage of the notions of cycles and currents. Let
$W$ be a scheme. We denote by $Z_k\ualg W$ the group of 
formal sums of finitely many
$k$-dimensional irreducible subvarieties of $W$ with rational
coefficients. We call elements of $Z_k\ualg W$ $k$-cycles of $W$. 
Now let $W$ be any stratified topological space with stratification $\cS$. 
We say that a (complex) $k$-dimensional current 
$C$ is stratifiable if there is a refinement of $\cS$, say $\cS\pri$, such that
there are finitely many $k$-dimensional strata $S_i$
and rationals $a_i\in\QQ$ such that
$C=\sum a_i[S_i]$ (All currents in this paper are oriented).
Here we assume that each stratum of
$\cS\pri$ was given an orientation a priori and $[S_i]$ is the oriented
current defined by $S_i$. We identify two currents if they define identical
measures in the sense of rectifiable currents.
We denote the set of all stratifiable
$k$-dimensional currents modulo the equivalence relation by $Z_k W$. 
Clearly, if $W$
is a scheme then any $k$-cycle has an associated current in
$Z_k W$, which defines a map $Z_k\ualg W\to Z_k W$.
In the following,  we will not
distinguish a cycle from its associated current.  
Hence for $C\in Z_k\ualg
W$ we will view it as an element of $Z_{k} W$. 
Note that if
$C\in Z_k  W$ has zero boundary in the sense of  current and $C$ has
compact support, then $C$ defines canonically an element in $H_{2k}(W,\QQ)$.
Finally, if $C=\sum a_i[S_i]\in Z_k W$ and $F\sub W$ is a stratifiable
subset, we say that $C$ intersects $F$ transversally if $F$ intersects each
$S_i$ transversally as stratified sets (See
\cite{GM} for topics on stratifications).
In such cases, we can define the intersection current $C\cap F$ if the
orientation of the intersection can be defined according to the geometry of
$W$ and $F$.     

We begin with a quick review of the algebraic construction of GW-invariants.
Let $X$ be any smooth projective variety and let $A\in H_2(X,\ZZ)$ and
$g,n\in\ZZ$ as before be fixed once and for all. We let $\mgna$ be the
moduli scheme of stable morphisms defined before. $\mgna$ is
projective. The GW-invariants of $X$ is defined using the
virtual moduli cycle 
$$[\mgna]\vir\in A\lsta\mgna.
$$
To review such a construction, a few
words on the obstruction theory of deformations of morphisms 
are in order. Let $w\in\mgna$ be any point associated to the
stable morphism $\cX$. Let $(B,I,\cX_{B/I})$ be any collection where $B$ is an
Artin ring, $I\sub B$ is an ideal annihilated by the maximal ideal 
$\mm_B$ of $B$ and $\cX_{B/I}$ is a flat family of stable morphisms
over $\spec B/I$ whose restriction to the closed fiber of $\cX_{B/I}$
is isomorphic to $\cX$. An obstruction theory to
deformation of $\cX$ consists of a $\CC$-vector space $V$, called the 
obstruction space, and an
assignment that assigns any data $(B,I,\cX_{B/I})$ as before to an obstruction
class 
$$\ob(B,B/I,\cX_{B/I})\in I\otimes_{\CC} V
$$
to extending $\cX_{B/I}$ to $\spec B$. Here by an obstruction class, 
we mean that its vanishing is the necessary and sufficient condition for
$\cX_{B/I}$ to be extendable to a family over $\spec B$.
We also require that such an assignment satisfies the obvious base change 
property (For reference on obstruction theory please consult 
\cite{Ob}). In case $\cX$ is the map
$f\mh C\to X$ with
$D\sub C$ the divisor of its $n$ marked points, the space of
the first order deformations of $\cX$ is parameterized by
$\Ext^1(\dbul_{\cX},\cO_C)$, where $\dbul_{\cX}=
[f\sta\om_X\to\om_C(D)]$ is the complex as before,
and the standard obstruction theory to deformation of $\cX$ 
takes values in $\Ext^2(\dbul_{\cX},\cO_C)$. 

An example of obstruction theories is the following.
Let $R$ be the ring of formal power series in $m$ variables
and let $\mm_R\sub R$ be its maximal ideal.
Let $F$ be a vector space and let $f\in \mm_R\otimes_{\CC}F$.
We let $(f)\sub R$ be the ideal generated by components of $f$. Then
there is a standard obstruction theory to deformations of
$0$ in $\spec R/(f)$ taking values in $V$, where $V$ is the 
cokernel of $df\mh (\mm_R/\mm_R^2)\dual\to F$,
defined as follows. Let $I\sub B$ be an ideal of an Artin ring as before and
let $\varphi_0\mh\spec B/I\to\spec R/(f)$
be any morphism. To extend $\varphi_0$ to $\spec B$, we first pick a
homomorphism $\sigma\mh R\to B$ extending the induced $R\to B/I$,
and hence a morphism $\varphi_{\text{pre}}\mh \spec B \to \spec R$. 
The image $\sigma(f)\in B\otimes F$ is in $I\otimes F$,
and is the obstruction to 
$\varphi_{\text{pre}}$ factor through $\spec R/(f)\sub\spec R$. 
Let $\ob(B,B/I,\varphi_0)$ be the image of $\sigma(f)$ in 
$I\otimes V$ via $F\to V$. It is direct to check that 
$\ob(B,B/I,\varphi_0)=0$ if and only if there is an extension 
$\varphi\mh\spec B\to\spec R/(f)$ of $\varphi_0$.
This assignment
\begin{equation}
\lab{eq:3.0}
(B,B/I,\varphi_0)\mapsto \ob(B,B/I,\varphi_0)\in
 I\otimes V
\end{equation}
is the induced obstruction theory of $\spec R/(f)$.

\begin{defi}
\lab{3.1}
A Kuranishi family of the standard obstruction theory of $\cX$ consists of
a vector space $F$, a ring of formal power series $R$ 
with $\mm_R$ its maximal ideal,
an $f\in\mm_R\otimes F$, a family $\cX_{R/(f)}$
of stable morphisms over $\spec R/(f)$ whose closed fiber over $0\in
\spec R/(f)$ is isomorphic to $\cX$ and an exact sequence
\begin{equation}
0\lra \Ext^1(\dbul_{\cX},\cO_C)\mapright{\alpha} (\mm_R/\mm_R^2)\dual
\mapright{df} F\lra \Ext^2(\dbul_{\cX},\cO_C)
\lra 0
\lab{eq:3.1}
\end{equation}
of which the following holds:
First, the composite
$$\Ext^1(\dbul_{\cX},\cO_C)\mapright{\alpha} \ker(df)\equiv T_0\spec 
R/(f) \mapright{} \Ext^1(\dbul_{\cX},\cO_C),
$$ 
where the second arrow is the Kodaira-Spencer map of the family $\cX_{R/(f)}$,
is the identity homomorphism;
Secondly, let $I\sub B$ and $\varphi_0\mh\spec B/I\to
\spec R/(f)$ be as before and let
$$\ob(B,B/I,\varphi_0\sta\cX_{R/(f)})\in I\otimes\Ext^2(\dbul_{\cX},\cO_C)
$$
be the obstruction to extending $\varphi_0\sta\cX_{B/I}$ to $\spec B$.
Then it is identical to $\ob(B,B/I,\varphi_0)$ 
under the isomorphism 
$$\coker(df)\cong \Ext^2(\dbul_{\cX},\cO_C),
$$
where $\ob(B,B/I,\varphi_0)$ is the obstruction class in \eqref{eq:3.0}.
\end{defi}

We now sketch how the virtual moduli cycle $[\mgna]\vir$ was constructed.
Similar to the situation of the moduli of stable smooth maps,
we need to treat $\mgna$ either as a $\QQ$-scheme or as a
Deligne-Mumford stack. The key ingredient here is the notion of atlas,
which is a collection of charts of $\mgna$. A chart of $\mgna$ 
is a tuple $(S, G, \cX_S)$, where $G$ is a finite group,
$S$ is a $G$-scheme 
(with effective $G$-action) and $\cX_S$ is a $G$-equivariant family of stable
morphisms so that the tautological morphism $\iota \mh S/G\to\mgna$ induced by
the family $\cX_S$ is an \'etale neighborhood.
For details of such an notion, please consult \cite{DM, Vi, LT1}.
We now let $f\mh C\to X$ be the representative of
$\cX_S$ with $D\sub C$ the divisor of the $n$-marked sections of $\cX_S$.
Let $\pi\mh C\to S$ be the projection. We consider the relative extension sheaves
$\ext^i_{\pi}(\dbul_{\cX_S},\cO_C)$, where 
$\dbul_{\cX_S}=[f\sta\om_X\to\om_{C/S}(D)]$ as
before. For short, we denote the sheaves 
$\ext^i_{\pi}(\dbul_{\cX_S},\cO_C)$ by $\cT^i_S$.
Because they vanish for $i=0$ and $i> 2$, for any $w\in S$, the
Zariski-tangent space $T_wS$ is $\cT_S^1\otimes_{\cO_S}k(w)$
and the 
obstruction space to deformations of $w$ in $S$ is $V_w=\cT_S^2\otimes
_{\cO_S}k(w)$. Now we choose a complex
of locally free sheaves of $\cO_S$-modules $\ebul=[\cE_1\to\cE_2]$ so that
it fits into the exact sequence
\begin{equation}
0\lra \cT_S^1\lra\cE_1\lra\cE_2\lra
\cT_S^2\lra 0.
\lab{eq:3.2}
\end{equation}
We let $F_i(w)=\cE_i\otimes_{\cO_S}k(w)$. 
Then we have the exact sequence of vector spaces
\begin{equation}
0\lra T_wS\lra F_1(w)\lra F_2(w) \lra V_w\lra 0.
\lab{3.31}
\end{equation}
We let $K_w\in  R(w)$ be a Kuranishi map of the
obstruction theory to deformations of $w$, where $R(w)=
\varprojlim\oplus _{k=0}^N \text{Sym}^k(F_1(w)\dual)$,
so that \eqref{3.31} is part of the data of the Kuranishi family specified
in Definition \ref{3.1}.
Let $(K_w)\sub R(w)$ be the ideal generated by the
components of $K_w$ and let $\spec R_w/(K_w)\sub\spec R_w$
be the corresponding subscheme. It follows that $\spec R_w/(K_w)$ is 
isomorphic to the formal completion of $S$ along $w$, denoted $\hat w$. 
We let $N_w$
be the normal cone to $\spec R_w/(K_w)$ in $\spec R_w$. Then $N_w$ is
canonically a subcone of $F_2(w)\times\hat w$. Here, by abuse of notation
we will use $F_2(w)$ to denote the total space of the vector 
space $F_2(w)$.
Note that $N_w$ is the infinitesimal normal cone to $S$ in its obstruction theory at $w$.
To obtain a global cone over $S$, we need the following existence 
and uniqueness theorem, which is the main result of [LT1].

In this paper, we will call a vector bundle $E$ the associated vector
bundle of a locally free sheaf $\cE$ if $\cO(E)\cong \cE$. For
notational simplicity, we will not distinguish a vector bundle from 
the total space (scheme)
of this vector bundle.

\begin{theo}[\cite{LT1}]
\lab{3.15}
Let $E$ be the associated vector bundle of $\cE_2$.
Then there is a cone scheme $N_S\sub E$ such that 
for each $w\in S$ there is an isomorphism
\begin{equation}
F_2(w)\times\hat w\cong E \times_S\hat w
\lab{eq:3.3}
\end{equation}
of cones over $\hat w$
extending the canonical isomorphism
$F_2(w) \cong E\times_S w$ such that under
the above isomorphism $N_w$ is 
isomorphic to $N_S\times_S\hat w$.
In particular, the cycle defined by the scheme $N_S$ is uniquely
characterized by this condition.
\end{theo}

In the previous discussion, if we replace $F_1(w)$ and $F_2(w)$ by
$T_wS$ and $V_w$ respectively, we obtain a Kuranishi map and
correspondingly a cone scheme in $V_w\times \hat w$, denoted
by $N_w^0$. 

\begin{theo}[\cite{LT1}]
\lab{3.16}
Let the notation be as before. Then there is a vector bundle homomorphism
$r\mh E\times_S\hat w\to V_w\times \hat w$ extending the canonical
homomorphism $E|_w\to V_w$ induced by \eqref{3.31}
such that 
$$N_w^0\times_{V_w\times\hat w} E\times_S\hat w= N_S\times_S \hat w.
$$
\end{theo}

To construct the virtual cycle $[\mgna]\vir$, we need to find a 
global complex over $\mgna$ analogous to $\ebul$. 
For the purpose of comparing with the analytic construction of
the virtual cycles, we will use atlas of analytic charts. 
We let $\{(R_i,W_i,\phi_i)\}_{i\in\Lambda}$ be the good
atlas of the smooth approximation of $[\Phi\mh\bB\to\bE]$ chosen 
in section 2.
Then the collection $Z_i=\phi_i\upmo(0)$ with the tautological family of 
stable analytic maps (with the last $k$-marked points discarded)
form an atlas of the underlying analytic scheme of $\mgna\cong\Phi\upmo(0)$.
Since we are only interested in constructing and working
with cone cycles in $\QQ$-bundles (known as V-bundles) over
$\mgna$, there is no loss of generality that we work with $\mgna$ with the
reduced scheme structure. Hence, for simplicity we will endow $Z_i=\phi_i\upmo(0)$ with the reduced analytic scheme structure. 
We let $\cX_i$ be the tautological
family of the $n$-pointed stable analytic maps over $Z_i$
that is derived by discarding the last $k$ marked points of
the restriction to $Z_i$ of the tautological family over $\tilde U_i$.
We let $G_i$ be the finite group associated to the
chart $(R_i,W_i,\phi_i)$, and let $\cX_i$ be represented by
$f_i\mh C_i\to X$ with $D_i\sub C_i$ be the divisor of the $n$-marked sections
of $C_i$ and let $\pi_i\mh C_i\to Z_i$ be the projection. 
In \cite{LT1}, to each $i$, 
we have constructed a $G_i$-equivariant complex
of locally free sheaves of $\cO_{Z_i}$-modules $\ebul_i=[\cE_{i,1}
\to\cE_{i,2}]$ such that $\ext_{\pi_i}^{\bullet}(
\dbul_{\cX_i},\cO_{C_i})$ is the sheaf cohomology of $\ebul_i$.
It follows from the algebraic and the analytic constructions of charts that
each $(Z_i,\cX_i)$ can be realized as an analytic open subset of
an algebraic chart, say $(S, G, \cX_S)$, and the complex $\ebul_i$ is the 
restriction to this open subset  of an algebraic complex $\ebul$, as in
\eqref{3.2}. Therefore we can apply Theorem \ref{3.15} to obtain a unique 
analytic cone cycle $M_i\ualg\in Z\lsta E_i$, where $E_i$ is the associated
vector bundle of $\cE_{i,2}$.
Let $\iota_i\mh Z_i/G_i\to\mgna $ be the tautological map induced by 
the family $\cX_i$. One property that follows from the construction of the
complexes $\ebul_i$ which we did not mention is that to each $i$, the 
cone bundle $E_i/G_i$ over $Z_i/G_i$ descends to a 
cone bundle over $\iota_i(Z_i/G_i)$, denoted by $\tilde E_i$, 
and 
$\{\tilde E_i\}_{i\in\Lambda}$ patch together to form a global cone
bundle over $\mgna$, denoted by $\tilde E$. Further, by the uniqueness 
of the cone cycles $M_i\ualg\in Z\lsta  E_i$ in
Theorem 3.2 and 3.3, to each $i$ the cone cycle $M_i\ualg/G_i$ in $E_i
/G_i$ descends to a cone cycle $\cM_i\ualg\in Z\lsta \tilde E_i$, and 
$\{\cM_i\ualg\}_{i\in\Lambda}$ patch together to form a cone cycle in $Z\lsta\tilde E$,
denoted by $\cM\ualg$.
It follows from \cite{LT1} that $\tilde E$ is an algebraic cone 
over $\mgna$ and $\cM\ualg$ is an algebraic cone cycle in $\tilde E$.
In the end, we let $\eta_E\mh \mgna\to \tilde E$ be the zero section
and let 
$$\eta_E\sta\mh \{\text{algebraic cycles in}\ Z\lsta\tilde E\}\lra H\lsta(\mgna;\QQ)
$$
be the Gysin homomorphism. Then the virtual moduli cycle is 
$$[\mgna]\vir=\eta_E\sta[\cM\ualg]\in A\lsta\mgna.
$$

There is an analogous way to construct the $GW$-invariants of algebraic
varieties using analytic method. We continue to use the notion developed
in section 1. Let $(R,W,\phi)$ be a smooth approximation of
$[\Phi\mh\bB\to\bE]$ constructed in Lemma \ref{1.20}. Then we can 
construct a cone current in the total space of $W$ as follows.
Let $\Gamma_{t\phi}$ be the graph of $t\phi$ in $W$ and 
let $N_{0/\phi}$ be the limit current $\lim_{t\to\infty}\Gamma_{t\phi}$,
when it exists. Clearly, if such a limit does exist, then it is
contained in $W|_{\phi\upmo(0)}$.
In general, though $\phi$ is smooth there is no guarantee
that such a limit will exist. However, if the approximation is analytic,
then we will show that such limit does exist as an
stratifiable current. Indeed, assume $(R,W,\phi)$ is an analytic smooth
approximation. Since the existence of $\lim\Gamma_{t\phi}$
is a local problem, we can
assume that there is a holomorphic basis of $W$, say $e_1,\ldots,e_r$. 
Then $\phi$ can be expressed
in terms of $r$ holomorphic functions $\phi_1,\ldots,\phi_r$.
Now let $\bC$ be the complex line with complex variable $t$, 
let $w_i$ be the dual of $e_i$ and let $\Theta\sub W\times \bC$ be the
analytic subscheme defined by the vanishing
of $tw_i-\phi_i$, $i=1,\ldots,r$.
We let $\Theta_0$ be the smallest closed analytic subscheme
of $\Theta$ that contains $\Theta\cap (W\times\bC\sta)$, where $\bC\sta
=\bC-\{0\}$. By the Weierstrass preparation theorem, such $\Theta_0$ does
exist. Then we define $N_{0/\phi}$ to be the associated cycle of
the intersection of the
scheme $\Theta_0$ with $W\times\{0\}$. By \cite{Fu},
$N_{0/\phi}$ is the limit of $\Gamma_{t\phi}$.
Obviously, $N_{0/\phi}$ is stratifiable.
This shows that for any analytic approximation $(R,W,\phi)$ the
limit $\lim\Gamma_{t\phi}$ does exist. 

We now state a simple lemma which implies that if $(R\pri,W\pri,
\phi\pri)$ is a smooth approximation that is finer than the
analytic approximation $(R,W,\phi)$, then $\lim\Gamma_{t\phi}$ exists as
well. We begin with the following situation. Let $V$ be a smooth
oriented vector bundle over a smooth oriented manifold $M$
and let $\varphi\mh M\to V$ be a smooth section. Let $V\pri\sub
V$ be a smooth submanifold such that for any $x\in\varphi\upmo(0)$
we have $\image(d\varphi(x))+V\pri_x=V_x$.
Then $M_0=\varphi\upmo(V\pri)$ is a smooth submanifold of
$M$ near $\varphi\upmo(0)$. Let $V_0$ be the restriction of
$V\pri$ to $M_0$ and let $\varphi_0\mh M_0\to V_0$ be
the induced section. We next let $N\sub TM|_{\varphi\upmo(0)}$ be a
subbundle complement to $TM_0|_{\varphi\upmo(0)}$
in $TM|_{\varphi\upmo(0)}$.
Then the union of $d\varphi(x)(N_x)$ for all $x\in\varphi\upmo(0)$
forms a subbundle of $V|_{\varphi\upmo(0)}$. We denote this bundle
by $d\varphi(N)$.
Since $V|_{\varphi\upmo(0)}\equiv V_0|_{\varphi\upmo(0)}
\oplus d\varphi(N)$, there is a unique projection
$P\mh V|_{\varphi\upmo(0)}\to V_0|_{\varphi\upmo(0)}$ such that
whose kernel is $d\varphi(N)$ and the composite
of the inclusion $V_0|_{\varphi\upmo(0)}\to V|_{\varphi\upmo(0)}$
with $P$ is the identity map. 

\begin{lemm}
\lab{3.2}
Let the notation be as before and let $l=\dim M$ and
$l_0=\dim M_0$. Then $\lim\Gamma_{t\varphi}$
exists as an $l$-dimensional current in $V|_{\varphi\upmo(0)}$ if and only if 
$\lim\Gamma_{t\varphi_0}$ exists as an $l_0$-dimensional 
oriented current in
$V_0|_{\varphi\upmo(0)}$. Further, if they do exist then
$$\lim\Gamma_{t\varphi}=P\sta(\lim\Gamma_{t\varphi_0}).
$$
Hence $\lim\Gamma_{t\phi}$ is stratifiable if $\lim\Gamma_{
t\phi_0}$ is stratifiable.
\end{lemm}

\begin{proof}
This is obvious and will be left to the readers.
\end{proof}

Now let $\{(R\lalp,W\lalp,\phi\lalp)\}_{\alpha\in\Xi}$ be a collection of
analytic smooth approximations of $[\Phi\mh\bB\to\bE]$ such
that the images of $Z\lalp=\phi\lalp\upmo(0)$ 
(in $\Phi\upmo(0)$) covers $\Phi\upmo(0)$.
It follows that we can choose a good atlas 
$\{(R_i,W_i,\phi_i)\}_{k\in\Lambda}$
constructed in Lemma \ref{1.20} so that all
approximations in $\Lambda$ are finer than approximations in $\Xi$.
Now let $i\in\Lambda$ and let $x\in Z_i=\phi_i\upmo(0)\sub R_i$ be any
point. Because charts in $\Xi$ cover $\Phi\upmo(0)$, there is an
$\alpha\in\Xi$ such that the image of $R\lalp$ in $\bB$ contains
the image of $x$ in $\bB$. Then because $(R_i,W_i,\phi_i)$ is finer than
$(R\lalp,W\lalp,\phi\lalp)$, by definition, there is a locally closed
submanifold $R_{i,\alpha}\sub R_i$, a local diffeomorphism $f_i^{\alpha}\mh
R_{i,\alpha}\to R\lalp$ and a vector bundle inclusion $(f_i^{\alpha})\sta W\lalp
\sub W_i|_{R_{i,\alpha}}$ such that $(f_i^{\alpha})\sta(\phi\lalp)=\phi_i$,
as in \eqref{eq:1.4}. 
This is exactly the situation studied in Lemma \ref{3.2}. Hence
$\lim\Gamma_{t\phi_i}$ exists near fibers of $W$ over $x$. 
Because $\{Z\lalp\}$ covers $\Phi\upmo(0)$,
$\lim\Gamma_{t\phi_i}$ exists and is a pure dimensional stratifiable
current of dimension $\dim R_i$. We denote this current by $N_{0/\phi_i}$.

Now it is clear how to construct the GW-invariants of algebraic
varieties using these analytically constructed cones. 
By the property of good coverings,
for $j\leq i\in\Lambda$ the approximation $(R_i,W_i,\phi_i)$ is finer than 
$(R_j,W_j,\phi_j)$.
We let $Z_i=\phi_i\upmo(0)$ be as before and let $Z_{i,j}=
Z_i\cap R_{i,j}\sub Z_i$, where $R_{i,j}$ is defined before \eqref{eq:1.4}.
Let $\rho^j_i\mh Z_{i,j}\to Z_j$ be the restriction of $f^j_i$ to $Z_{i,j}$.
Note that $Z_{i,j}$ is an open subset of $Z_i$ and $\rho^j_i\mh Z_{i,j}\to 
\rho^j_i(Z_{i,j})$ is a covering. 
Let $F_i$ be the restriction of $W_i$ to $Z_i$ and let $p_i\mh F_i\to Z_i$
be the projection.  Hence, $(\rho_i^j)\sta(F_j)$ is canonically
a subvector bundle of $F_i|_{Z_{i,j}}$. 
By Lemma \ref{3.2}, $(\rho^j_i)\sta(F_j)$ intersects
transversally with $N_{0/\phi_i}\cap p_i\upmo(Z_{i,j})$ and as currents,
$N_{0/\phi_i}\cap (\rho_i^j)\sta(F_j)=(p_i^j)\sta(N_{0/\phi_j})$.

For convenience, in the following we will call the collection 
$\{F_i\}$ with transition functions $f^j_i$ a semi-$\QQ$-bundle and
denote it by $\cF$, and will denote $\{ N_{0/\phi_i}\}$ by $\cN\uan$.
As in section two, we call a collection 
$\bs=\{s_i\}_{i\in\Lambda}$ of smooth sections $s_i\mh Z_i\to F_i$
a global section of 
$\cF$ if for $j\leq i\in\Lambda$ the restriction $s_i|_{Z_{i,j}}
\mh Z_{i,j}\to F_i|_{Z_{i,j}}$ coincides with the pull back section
$(\rho_i^j)\sta s_j\mh Z_{i,j}\to (\rho_i^j)\sta F_j$ under the canonical 
inclusion $(\rho_i^j)\sta F_j\sub F_i|_{Z_{i,j}}$. We say that the section 
$\bs$ is transversal to  $\cN\uan$ if for each $i\in\Lambda$,
the graph of the section $s_i$ is transversal to $N_{0/\phi_i}$ in $F_i$.

Obviously, if $\bs$ is a global section of $F$
that is transversal to $\cN\uan$, then following the argument
after Lemma \ref{1.10}, currents 
$$\frac{1}{m_i}\iota\pri_{i\ast}\pi_{i\ast}(N_{0/\phi_i}\cap
\Gamma_{s_i}),\quad {i\in\Lambda},
$$
 where $\iota\pri_i\mh Z_i\to\bB$ is the restriction of 
$\iota_i\mh R_i\to\bB$ to $Z_i\sub R_i$ and $m_i$ is the number of sheets 
of the branched covering $\iota\pri_i\mh Z_i\to\iota_i\pri(Z_i)$, patch 
together to form an oriented current in $\bB$ without boundary.
We denote this current by $\bs\sta(\cN\uan)$. 
It has pure dimension $r\lexp$ since the currents
$N_{0/\phi_i}$ has dimension $\dim R_i=\rank F_i+r\lexp$.
Hence it defines a homology class $[\bs\sta(\cN\uan)]$
in $H_{2r\lexp}(\bB;\QQ)$.

\begin{prop}
$[\bs\sta(\cN\uan)]$ is the Euler class $e[\Phi\mh \bB\to\bE]$
constructed in section one.
\end{prop}

\begin{proof}
Recall that the class $e[\Phi\mh\bB\to\bE]$ was constructed
by first selecting a collection of perturbations $h_i(s)\mh R_i\to W_i$
of $\phi_i$ parameterized by $s\in [0,1]$
satisfying certain property and then form the current that is the
patch together of the currents $\frac{1}{ m_i}\iota_{k\ast}(\Gamma_{h_i(1)}
\cap\Gamma_0)$, where $\Gamma_{h_i(1)}$ and $\Gamma_0$ 
are the graph of $h_i(1)$ and $0\mh R_i\to W_i$. 
Alternatively, we can perturb the $0$-section instead of 
$\{\phi_i\}$ to obtain the same cycle. Namely, we let $h_i\pri(s)\mh
R_i\to W_i$ be a collection of perturbations of the zero section
$0\mh R_i\to W_i$, such that it satisfies the obvious
compatibility and properness property similar to that of
$h_i(s)$ in section two.
Moreover, we require that the graph $\Gamma_{h\pri_i(1)}$
is transversal to $N_{0/\phi_i}$ and transversal to
the graph $\Gamma_{t\phi_i}$ for sufficiently large $t$.
Of course such perturbations do exist following the proof of Proposition
\ref{1.12}. Let $C_t$ be the current in $\bB$ that is the result of
patching together the currents
$\frac{1}{ m_i}\iota_{i\ast}p_{i\ast}(\Gamma_{h_i\pri(1)}
\cap\Gamma_{t\phi_i})$, where $p_i$ is the projection $W_i\to R_i$. 
Clearly, for $t\gg 0$, we have $\partial C_t=0$ and 
$\supp(C_t)$ is compact. Hence $C_t$ defines a homology
class in $H_{2r\lexp}(\bB;\QQ)$, denoted by $[C_t]$.
It follows from the uniqueness argument in the end of section two
that for sufficiently large $t$, the homology class
$[C_t]$ in
$H_{2r\lexp}(\bB;\QQ)$ is exactly the Euler class.
On the other hand, we let $C_{\infty}$ be the current in
$\bB$ that is the patch together of the currents
$\frac{1}{ m_i}\iota_{i\ast}p_{i\ast}(\Gamma_{h_i\pri(1)}\cap N_{0/\phi_0})$.
Because $N_{0/\phi_i}$ is the limit of 
$\Gamma_{t\phi_i}$, and because $\Gamma_{h_i\pri(1)}$ intersects
transversally with $\Gamma_{t\phi_i}$ for $t\gg 0$ and with $N_{0/\phi_i}$,
the union
$$\bigcup_{t\in [0,\eps]} \{t\}\times C_{1/t}\sub [0,\eps]\times\bB,
$$
where $1\gg\eps>0$, is a current whose boundary is $C_{1/\eps}- C_{\infty}$.
This implies that
$$[C_{\infty}]=[C_t]\in H_{2r\lexp}(\bB;\QQ)\qquad \text{for}\ t\gg 0.
$$
Further, because the currents $N_{0/\phi_i}$
are contained in $F_i=W_i|_{Z_i}$,
$p_{i\ast}(N_{0/\phi_i}\cap
\Gamma_{h_i\pri(1)})$ as current is identical to 
$\pi_{i\ast}(N_{0/\phi_i}\cap \Gamma_{r_i})$, where $r_i\mh
Z_i\to F_i$ is the restriction of $h_i\pri(1)$ to $Z_i$. 
Hence $C_{\infty}={\mathbf r}\sta(\cN\uan)$ with ${\mathbf r}=\{r_i\}$.
Finally, it is direct to check that the homology classes
$[\bs\sta(\cN\uan)]$ do not depend on the choices of the 
sections $\bs$ of $\cF=\{F_i\}$ so long as they
satisfy the obvious transversality conditions.
Therefore, 
$$[\bs\sta(\cN\uan)]=[{\mathbf r}\sta(\cN\uan)]=
[C_{1/\eps}]=e[\Phi:\bB\to\bE].
$$
This proves the Proposition.
\end{proof}

In the end, we will compare the algebraic normal cones with the analytic
normal cones to demonstrate that the algebraic and analytic construction
of the GW-invariants give rise to the identical invariants.

Here is our strategy. Taking the good atlas
$\{(Z_i,\cX_i)\}_{i\in\Lambda}$ of $\mgna$ as before, we have two collections of
semi-$\QQ$-vector bundles, namely $\cE=\{E_i\}$ and $\cF=\{F_i\}$, and two
collections  of cone currents $\cM\ualg=\{M_i\ualg\}$ and
$\cN\uan=\{N_{0/\phi_i}\}$ such that 
$[\eta_E\sta(\cM\ualg)]$ and $[\eta_F\sta(\cN\uan)]$ are the algebraic and
the symplectic virtual moduli cycles of $\mgna$ respectively. 
Here $\eta_E$ and $\eta_F$ are generic sections of $\cE$ and $\cF$
respectively. To
compare these two classes, we will form a new 
semi-$\QQ$-vector bundle $\cV=\{V_i\}$, where
$V_i=E_i\oplus F_i$, and construct
a stratifiable cone current $\cP$ in $\cV$
such that the cycle $\cP$ intersect $\cE\sub\cV$ and 
$\cF\sub\cV$ transversally and the intersection $\cP\cap \cE$ and $\cP\cap
\cF$ are $\cM\ualg$ and $\cN\uan$ respectively. Therefore, if we let
$\eta_V$ be a generic section of $\cV$, then
$$[\eta_E\sta(\cM\ualg)]=[\eta\sta_V(\cP)]=[\eta\sta_F(\cN\uan)]
\in H\lsta(\mgna;\QQ).
$$
This will prove the Comparison Theorem.

We now provide the details of this argument. We begin with any index
$i\in\Lambda$ and an open subset $S\sub Z_i$.
Let $f\mh C\to X$ be the restriction to $S$ of the
tautological family $\cX_i$ of stable maps over $Z_i$, with
$D\sub C$ the divisor of its $n$-marked sections
and $\pi\mh C\to S$ the projection. 
Note that $f$ is the restriction of a family of stable 
morphisms over a scheme to an analytic open subset of the base scheme.
Following the construction in \cite[section 3]{LT1}, after fixing a
sufficiently ample line bundle over $X$, we canonically construct
a locally free sheaf of $\cO_{C}$-modules $\cK$ so that
$f\sta\om_X$ is canonically a quotient sheaf of $\cK$. Let
$\cL$ be the kernel of $\cK\to f\sta\om_X$. Then 
the restriction to $S$ of the sheaf
$\cE_{i,1}$ (resp. $\cE_{i,2}$)
mentioned before is the the relative extension
sheaf $\ext_{\pi}^1([\cK\to\om_{C/S}(D)],\cO_{C})$
(resp. $R\pi_{\ast}(\cL\dual)$). 
We denote them by $\cE_{S,1}$ and $\cE_{S,2}$ respectively.
As usual, we let $E_{S,1}$ and $E_{S,2}$ be the
associated vector bundle of $\cE_{S,1}$ and $\cE_{S,2}$
respectively.
Following the notation in \cite{LT1}, the
tangent-obstruction complex $[\cT_S^1\to\cT_S^2]$ of $\cX_i|_S$ is
$$\bigl[\sideset{}{^1_{\pi}}\ext([f\sta\om_X\to\om_{C/S}(D)],\cO_{C})
\mapto{\times0}
\sideset{}{^2_{\pi}}\ext([f\sta\om_X\to\om_{C/S}(D)],\cO_{C})\bigr],
$$
and that there is a canonical homomorphism
$\eps: \cE_{S,1}\lra \cE_{S,2}$
so that the kernel and the cokernel of $\eps$ are
$\cT_S^1$ and $\cT_S^2$ respectively.
The homomorphism $\eps$ is the middle arrow in the sequence 
\eqref{eq:3.2}. 

We now assume that there is an analytic approximation
$\alpha\in\Xi$ so that $(R_i,W_i,\phi_i)$ is finer than $\alpha$
and $\iota_i(S)\sub\bB$
is contained in $\iota\lalp(Z\lalp)$. 
Let $\rho_{\alpha}\mh Z_i\to Z\lalp$
be induced by $f_i^{\alpha}\mh R_{i,\alpha}\to R\lalp$ (see \eqref{eq:1.4}).
Let $F_{S,\alpha}$
be the vector bundle over $Z_i$ that is the pull back of $F\lalp$.
Note that $F_{S,\alpha}$ is a smooth vector bundle. Let 
$G_{S,\alpha,2}= E_{S,2}\oplus F_{S,\alpha}$.
In the following, we will construct a holomorphic
vector bundle $G_{S,\alpha,1}$ and a
possibly degenerate
vector bundle homomorphism $\beta$ and non-degenerate
vector bundle inclusions $\tau_{\alpha,j}$ as shown below so that
\begin{equation}
\begin{CD}
E_{S,1} @>{\eps}>> E_{S,2}\\
@VV{\tau_{\alpha,1}}V @VV{\tau_{\alpha,2}}V\\
G_{S,\alpha,1} @>{\beta}>>G_{S,\alpha,2}
\end{CD}
\lab{eq:3.5}
\end{equation}
is commutative.
Let $w$ be any point in $S$. We denote by $C_w$ 
the fiber of $C$ over $w$ and let $f_w$ (resp. $\cK_w$, resp.
$\cL_w$) be the restriction of the respective objects to $C_w$.
As before, for any locally free
sheaf of $\cO_{C_w}$-modules $\cW$ that is locally free
away from the nodal points of $C_w$, we denote by $\cW\ua$ the sheaf
whose stalk at nodal points $z$ of $C_w$ are $\cW_z$ and its stalks 
at smooth points $z$ of $C_w$ are germs of smooth sections of the
associated vector bundle of $\cW$ at $z$. We let $G_{S,\alpha,1}|_w$ 
be the vector space of the equivalence classes of
commutative diagrams
\begin{equation}
\begin{CD}
@. @. \cK_w @>>> f_w\sta\om_X\\
@.@. @VV{h}V @VV{df_w\dual}V\\
0@>>> \cO_{C_w}\ua @>>> \cB_w\ua @>>> \om_{C_w}(D_w)\ua
@>>> 0
\end{CD}
\lab{eq:3.6}
\end{equation}
such that the lower exact sequences are induced by the exact sequences 
of sheaves of $\cO_{C_w}$-modules
\begin{equation*}
\begin{CD}
0 @>>> \cO_{C_w} @>>> \cB_w @>>> \om_{C_w}(D_w)
@>>> 0
\end{CD}
\end{equation*}
and that $h$ satisfies the following two requirements. 
First, let 
$c\mh \cL_w\to \cB_w\ua$ be the composite of $\cL_w\to\cK_w$ with
$h$. Since $\cL_w$ is the kernel of $\cK_w\to f_w
\sta\om_X$, $c$ automatically lifts to $h_{E}\mh \cL_w
\to \cO_{C_w}\ua$. The first requirement is that $h_{E}$ is
holomorphic. Secondly, since both $\cK_w$ and $\cL_w$ are sheaves of 
$\cO_{C_w}$-modules and since $h$ is analytic near nodal points of $C_w$, 
$\dbar h$ is a $(0,1)$-form with compact support 
\footnote{By which we mean that $\dbar h$ vanishes in
a neighborhood of the nodal points of $C_w$.}
taking values
in the associated vector bundle of $\cK_w\dual\otimes_{\cO
_{C_w}}\cB_w$. Because of the first requirement, it factors through
a section $h_{F}$ of $\om^{0,1}\lcpt(f_w\sta T_X)$.
We require that $h_{ F}$ is an element in $\rho_{\alpha}\sta W\lalp|_w$.
Using Lemma \ref{2.1} and Corollary \ref{2.2} and
the fact that $\cK\dual$ is sufficiently ample which 
was the precondition on our choice of $\cK$,
it is direct to check that the collection
$\{G_{S,\alpha,1}|_w\mid w\in S\}$ forms a
smooth vector bundle, denoted $G_{S,\alpha,1}$, and the correspondence
that sends \eqref{eq:3.6} to $h_{E}- h_{F}$ 
form a possibly degenerate
vector bundle homomorphism
$\beta \mh G_{S,\alpha,1}\to G_{S,\alpha,2}$. 

We next define the
homomorphisms $\tau_{\alpha,j}$. The 
homomorphism $\tau_{\alpha,2}\mh E_{S,2}\to G_{S,\alpha,2}$
is the obvious homomorphism based on the definition
 $G_{S,\alpha,2}=E_{S,2}\oplus F_{S,\alpha}$. For
$\tau_{\alpha,1}$, we recall that for any $w\in S$ the vector space
$E_{S,1}|_w$ is the set of equivalence classes of the diagrams
\eqref{eq:3.6} of which the $h$ are holomorphic. Namely,
$h$ are induced by   homomorphisms
$f_w\sta\om_X\to\cB$. 
Hence $E_{S,1}$ is canonically a subbundle of $G_{S,\alpha,1}$. 
This shows that both
$\tau_{\alpha,1}$ and $\tau_{\alpha,2}$ are 
inclusions of vector bundles.
Finally, let $\xi\in E_{S,1}|_w$ be any element associated to the diagram
\eqref{eq:3.6}, then
$\eps(\xi)$ is the section of $\cL_w\dual$
that is the lift of $\cL_w\to\cK_w\mapto{h}\cB_w$ to $\cL_w\to\cO_{C_w}$. It
follows that the square of \eqref{eq:3.5} is commutative.
We now show that $\coker(\tau_{\alpha,1})=
\coker(\tau_{\alpha,2})$. It suffices to show that  the sequence
\begin{equation}
\begin{CD}
0 @>>>  E_{S,1}@>{\tau_{\alpha,1}}>> G_{S,\alpha,1} @>{c}>> F_{S,\alpha} @>>> 0,
\end{CD}
\lab{eq:3.9}
\end{equation}
where $c$ is the composite of $\beta$ with $G_{S,\alpha,2}\to F_{S,\alpha}$,
is an 
exact sequence. But this follows directly from the definition of 
$G_{S,\alpha,1}$ and Lemma \ref{2.1} and Corollary \ref{2.2}.
This proves that $\coker(\tau_{\alpha,1})=
\coker(\tau_{\alpha,2})$, and consequently
\begin{equation}
\coker(\beta|_w)=\coker(\eps|_w)=\cT^2_S|_w
\lab{eq:3.31}
\end{equation}
for any $w\in S$.

In the following, we will construct the cone current $Q\lsa\in Z\lsta
V\lsat$. We first pick a subbundle $H\lalp\sub G\lsat$ such that
$H\lalp\to G\lsat\to\coker(\tau_{\alpha,1})$ is an isomorphism. We let
$P\lalp\mh G_{S,\alpha,2}\to E_{S,2}$ be the projection so that
$\ker(P\lalp)=\beta(H\lalp)$ and
$P\lalp\circ \tau_{\alpha,2}={\mathbf
1}_{E\lst}$. We then take
$Q\lsa$ to be the flat pull back current $P\lalp\sta(M_i\ualg)\in
Z\lsta G_{S,\alpha,2}$. It follows that $Q\lsa$ 
intersects the subbundle $E_{S,2}
\sub G_{S,\alpha,2}$ transversally and the intersection
$Q\lsa\cap E_{S,2}$ is exactly $M\ualg_S=M_i\ualg|_S$.
In the following, we will demonstrate that $Q\lsa$ intersects
the subbundle $F_{S,\alpha}\sub G_{S,\alpha,2}$ transversally as well
and that the
intersection $Q\lsa\cap F\lsa$ is the current
$\rho\lalp\sta (N\lalp\uan)\in Z\lsta F_{S,\alpha}$. 
Let $w\in S$ ($\sub Z_i$) be any point. Since $T_{w\pri}R\lalp$, where
$w\pri=\rho_{\alpha}(w)$, is the vector space 
$\Ext^1(\dbul_w,\cO_{C_w}\ua)\shar$,
there is a canonical injective homomorphism 
$\sigma_w \mh T_{w\pri}R\lalp\to G_{S,\alpha,1}|_w$ of vector spaces
that send the diagram \eqref{eq:2.3} to \eqref{eq:3.6} with
$\cK_w\to \cB_w\ua$ the composite of
$\cK|_w\to f_w\sta\om_X$ with the $v_2$ in \eqref{eq:2.3}. 
It is easy to see that the collection $\{\sigma_w \}_{w\in S}$ 
forms a smooth
non-degenerate vector bundle homomorphism 
$\sigma \mh \rho_{\alpha}\sta(TR\lalp)\to G_{S,\alpha,1}$. If follows
from the description of 
$$\rho\lalp\sta(d \phi\lalp): 
\rho\lalp\sta(TR\lalp)\lra F_{S,\alpha}
$$
that the diagram of vector bundle homomorphisms
\begin{equation}
\begin{CD}
G_{S,\alpha,1} @>{\beta}>> G_{S,\alpha,2}\\
@AAA @AAA\\
\rho\lalp\sta(TR\lalp) 
@>{\rho\lalp\sta(d\phi\lalp)}>> F_{S,\alpha}
\end{CD}
\end{equation}
is commutative,
where the second vertical arrow is the obvious 
inclusion.

To compare $Q\lsa$ with $\rho\lalp\sta(N\lalp\uan)$, we need
the following two lemmas.

\begin{lemm}
\lab{3.10}
Let $w\in S$ be any point and let $w\pri=\rho_{\alpha}(w)$.
Let $d_2\mh G_{S,\alpha,2}|_w\to \cT_S^2|_w$
be the homomorphism induced by \eqref{eq:3.31} and let
$F_{S,\alpha}|_w\to
\cT_S^2|_w$ be the canonical homomorphism given in
Lemma \ref{2.22}.
Then the following squares are commutative:
\begin{equation}
\begin{CD}
F_{S,\alpha}|_w @>{\sub}>> G_{S,\alpha,2}|_w @<{\tau_{\alpha,2}}<<
E\lst\\
@V{d_3}VV @V{d_2}VV @V{d_1}VV\\
\cT_S^2|_w @= \cT_S^2|_w @=\cT_S^2|_w 
\end{CD}
\lab{eq:3.40}
\end{equation}
\end{lemm}

\begin{lemm}
\lab{3.11}
For any point $w\in Z\lalp$, the germ of $\phi\lalp\mh R\lalp\to
W\lalp$ at $w$ is a Kuranishi map of the standard obstruction theory
of the deformation of stable morphisms associated to the exact sequence
$$
0\lra \cT\lalp^1|_w\lra T_wR\lalp \lra F\lalp|_w\lra
\cT\lalp^2|_w\lra 0.
$$
\end{lemm}

\begin{proof}
We first prove Lemma \ref{3.10}.
Since $G\lsat\equiv E\lst\oplus F\lsa$, $d_1$ and $d_3$ induces
a homomorphism $G\lsat|_w\to \cT^2_S|_w$. To prove
the lemma, it suffices to show that $d_2=d_1\oplus d_3$.
To accomplish this, we only need to show that for any 
$\xi\in G_{S,\alpha,1}|_w$ with
$\xi_E$ and $-\xi_F$ its two components of
$\beta(\xi)$ according to the direct sum
decomposition $G\lsat|_w=E_{S,2}|_w\oplus F_{S,\alpha}|_w$,
then
$d_1(\xi_E)=d_3(\xi_F)$.

To prove this, we first pick an $h_0\mh f_w\sta\om_X\to \cB\ua_w$ such that
\begin{equation}
\begin{CD}
f_w\sta\om_X @= f_w\sta\om_X\\
@VV{h_0}V @VV{df_w\dual}V\\\
\cB\ua_w @>>> \om_{C_w}(D_w)\ua
\end{CD}
\end{equation}
is commutative. Let $h_0\pri$ be the composite of $\cK_w\to f_w\sta\om_X$
with $h_0$.
Then $h\pri-h_0$ factor through $\cO_{C_w}\ua\to\cB_w\ua$,
say $\tilde h\mh \cK_w\to\cO_{C_w}\ua$. Clearly,
$\tilde h$ composed with $\cL_w\to\cK_w$ is the section
$\xi_E\in H^0(\cL_w\dual)$.
On the other hand, the lift of $\dbar \tilde h$ to 
$\om^{0,1}\lcpt(f_w\sta T_X)$
is $\xi_F-(\dbar h_0)^{\text{lift}}$. By the definition of
the connecting homomorphism $\delta\mh H^0(\cL_w\dual)\to
H^1(f_w\sta\om_X\dual)$, 
$$\delta(\xi_E)=\text{the image of}\ (\xi_F
-(\dbar h_0)^{\text{lift}})\ \text{in}\
H^{0,1}_{\dbar}(f_w\sta T_X)\cong H^1(f_w\sta\om_X\dual).
$$
However, the image of $(\dbar h_0)^{\text{lift}}$ is contained in the image of
the connecting homomorphism
$$ \Ext^1(\om_{C_w}(D_w),\cO_{C_w})\lra 
\Ext^2([f_w\sta\om_X\to 0],\cO_{C_w})\equiv H^1(f_w\sta\om_X\dual).
$$
Hence $d_1(\xi_E)=d_3(\xi_F)$. This proves Lemma \ref{3.10}.
\end{proof}

\begin{proof}
We now prove Lemma \ref{3.11}.
Let $I\sub B$ be an ideal of an Artin ring annihilated by the
maximal ideal $\mm_B$ and let $\varphi\mh\spec B/I
\to R\lalp$ be a morphism that sends the closed point of 
$\spec B/I$ to $w$ and
such that $\varphi\sta(\phi\lalp)=0$. By the description of 
the tautological family $\cX\lalp$ over $R\lalp$, the
pull back $\varphi\sta(\cX\lalp)$ forms an algebraic family 
of stable morphisms over
$\spec B/I$. We continue to use the open covering of the domain
$\cX\lalp$ used before. Since $R\lalp$ is smooth, we can extend
$\varphi$ to $\tilde\varphi\mh\spec B\to R\lalp$. Let $C_B$ 
over $\spec B$ be the domain of the pull back of the domain of
$\cX\lalp$ via $\tilde\varphi$ and
let $C_{B/I}$ be the domain of $C_B$ over $\spec B/I$. We let 
$\{U_i\}$ (resp. $\{\tilde U_i\}$)
be the induced open covering of $C_{B/I}$ (resp. $C_B$)
and let $f_i\mh U_i\to X$
be the restriction to $U_i$ of the pull back of the stable maps in
$\cX\lalp$. Because $\varphi\sta(\phi\lalp)=0$, $f_i$ are 
holomorphic. Hence they define a morphism $f\mh C_{B/I}\to X$.
Now we describe the obstruction to extending $f$ to $\spec B$. Let
$C_0$ be the closed fiber of $C_B$ and let $f_0\mh C_0\to X$
be the restriction of $f$. For each $i$, we pick a holomorphic
extension $\tilde f_i\mh \tilde U_i\to X$ of $f_i$. Then over 
$\tilde U_{ij}=\tilde U_i\cap \tilde U_j$, $\tilde f_j-\tilde f_i$
is canonically an element in $\Gamma(f_0\sta\cT_X|_{U_i\cap U_j})\otimes I$, 
denoted by
$f_{ij}$. Further, the collection $\{f_{ij}\}$ 
is a cocycle and hence defines an 
element $[f_{ij}]\in H^1(f_0\sta\cT_X)\otimes I$.
The obstruction to extending $f$ to $\spec B$ is the image of
$[f_{ij}]$ in $\Ext^2(\cD^{\bullet}_w,\cO_{C_0})\otimes I$
under the homomorphism in the statement in Lemma \ref{2.22}
with $z$ replaced by $w$.
We denote the image by $\text{ob}\ualg$.

The obstruction to extending $\varphi$ to 
$\tilde\varphi\mh\spec B\to
R\lalp$ so that $\tilde\varphi\sta(\phi\lalp)=0$ can be constructed
as follows. Let $g_i\mh \tilde U_i\to X$ be the pull back of 
the maps in $\cX\lalp$.
Note that $g_i$ are well defined since maps in $\cX\lalp$ depend 
analytically on the base manifold $R\lalp$. By the construction of
$R\lalp$, for each $i>l$ the map $g_i$ is
holomorphic. For $i<l$, we have canonical biholomorphism 
$\tilde U_i\cong \spec B\times (U_i\cap C_0)$. 
Because $\varphi\sta(\phi\lalp)\equiv0$, if we 
let $\xi_i$ be a holomorphic variable of $U_i\cap C_0$, then
$\frac{\partial}{\partial\bar\xi_i}g_i\cdot d\bar\xi_i$,
denoted in short $\dbar g_i$, vanishes
over $U_i\sub \tilde U_i$. Hence $\dbar h$ is a section of 
$\Gamma(\om^{0,1}\lcpt(f_0\sta T_X)|_{U_i\cap C_0})\otimes I$.
Clearly they patch together to form a global section
$\gamma$ of $\om\lcpt^{0,1}(f_0\sta T_X)\otimes I$.
The element $\gamma$ can be also defined as follows. 
Let ${\tilde\varphi}\sta \mh \cO_{R\lalp}\to B$ be the induced homomorphism
on rings. Then since the image of ${\tilde\varphi\lalp}\sta(\phi\lalp)\in
B\otimes_{\cO_{R\lalp}}\!\cO_{R\lalp}(W\lalp)$ in
$B/I\otimes_{\cO_{R\lalp}}\!\cO_{R\lalp}(W\lalp)$ vanishes, 
it induces an element
$\gamma\pri\in I\otimes W\lalp|_w$.
By our construction of $R\lalp$ and $\phi\lalp$, $\gamma$
coincides with $\gamma\pri$
under the inclusion $W\lalp|_w\sub 
\Gamma_{C_0}(\om\lcpt^{0,1}(f_0\sta T_X))$.
Let $\text{ob}\uan$ be the
image of $\gamma$ in the cokernel of 
$d\phi\lalp(w)\mh T_w R\lalp\to W\lalp|_w$. By definition,
$\text{ob}\uan$ is the obstruction to extending
$\varphi$ to $\tilde\varphi\mh \spec B\to \{\phi\lalp=0\}$.

To finish the proof of the lemma, we need to show that 
$\text{ab}\ualg=\text{ob}\uan$ under the isomorphism
$$\coker\{d\phi\lalp(w)\}\cong \Ext^1(\dbul_w,\cO_{C_0})
$$
given in Lemma \ref{2.22}. For this, it suffices to show that the Dolbeault
cohomology class of $\gamma$, denoted $[\gamma]\in H^{0,1}_{\dbar}(f_0\sta
T_X)\otimes I$, coincides with the C\v{e}ch cohomology class $[f_{ij}]\in
H^1(f\sta_0\cT_X)
\otimes I$ under the canonical isomorphism $H^{0,1}_{\dbar}(f_0\sta
T_X)\cong H^1(f\sta_0\cT_X)$. But this is obvious since $\varphi_i
=\tilde f_i-g_i$ is in $\Gamma_{U_i\cap C_0}(
\om\lcpt^0(f_0\sta T_X))\otimes I$ such that 
$\varphi_j-\varphi_i=f_{ij}$
and $\dbar \varphi_i=-\dbar g_i$.
Hence, $[f_{ij}]=[\gamma]$ under the given 
isomorphism. This proves the lemma.
\end{proof}

Now we come back to $Q\lsa\in Z\lsta G\lsat$. Let $w\in S$
be any point, let $\hat w$ be the formal completion of $S$ along $w$, let $V_w$
be the total space of $\cT^2_S|_w$ and let $N_w^0\sub V_w\times \hat w$ be the 
the cone in Theorem \ref{3.16}. We let $M_S\ualg$, 
$N\lsa\uan=\rho\lalp\sta(N\uan_i)$ and $Q\lsa$ be the cone 
currents in $E\lst$, $F\lsa$ and $G\lsat$ respectively as before.
Note that they are supported on union of closed subsets each diffeomorphic to
analytic variety. By Theorem \ref{3.15}, we have
vector bundle homomorphisms
$$e_1: E\lst\times_S\hat w\lra V_w\times \hat w
\and
e_3: F\lsa\times_S\hat w\lra V_w\times \hat w
$$
extending $E\lst|_w\to\cT^2_S|_w$ and $F\lsa|_w\to\cT^2_S|_w$ such that
$e_1\sta(N_w^0)$ and $e_3\sta(N_w^0)$ are the restrictions of $M_S\ualg$
and $N\lst\uan$ to fibers over $\hat w$ in $S$ respectively. Let 
$e_2\mh G\lsat\times_S\hat w
\to V_w\times\hat w$ be induced by $P\lalp \mh G\lsat\to E\lst$ and $e_1$.
Then $e_2\sta(N_w^0)$ is the restriction of $Q\lsa$ to $G\lsat\times_S\hat w$.
Because the squares in \eqref{eq:3.40} are commutative, 
\begin{equation*}
\begin{CD}
e_2\mh F\lsa\times_S\hat w
@>{\sub}>> G_{S,\alpha,2}\times_S\hat 
w @>{e_2|_{\hat w}}>> V_w\times\hat w
\end{CD}
\end{equation*}
is surjective. Hence $F\lsa\times_S\hat w$
intersects $Q\lsa$ transversally along fiber over $w$.
Let $e_3\pri\mh F\lsa\times_S\hat w\to V_w\times\hat w$ be induced by 
$F\lsa\to G\lsat$ and $e_2$, then the intersection of $Q\lsa$ with 
$F\lsa\times_S\hat w$ is $(e_3\pri)\sta(N_w^0)$.
However, by the choice of $P\lalp$, we have
$e_3\pri\equiv e_3|_w$, therefore the support
of $Q\lsa\cap F\lsa|_w$ is identical to the support of $N\lsa\uan|_w$.
Because $w\in S$ is arbitrary, the support of $Q\lsa\cap F\lsa$ is
identical to the support of $N\lsa\uan$. 
Further, for the same reason, for general point $p$ in $N\lsa\uan$
the multiplicity of $N\lsa\uan$ at $p$ is identical to the
multiplicity of the corresponding point in $Q\lsa\cap F\lsa$. 
This proves that
the cycles (or currents) $Q\lsa$ intersect $F\lsa\sub G\lsat$ transversally
and $Q\lsa\cap F\lsa=N\lsa\uan$. We remark that for the same reason, the
current $Q\lsa$ is independent of the choice of the subbundles $H\lalp\sub
G\lsat$.

We now let $F_S=F_i|_S$ and let $G\lst=E\lst\oplus F_S$. Note that
$G\lsat\sub G\lst$.
Because $R_i$ is finer than $R\lalp$, $\rho\lalp\sta TR\lalp$ is
a subbundle of $TR_i|_S$. Let $K\lalp\sub TR_i|_S$ be a complement of 
$\rho\lalp\sta TR\lalp\sub TR_i|_S$ and let $d\phi_i(K\lalp)\sub F_S$
be the image of this subbundle. Let $P_{S,\alpha}\mh F_S\to F\lsa$ be the
projection so that $\ker P_{S,\alpha}=d\phi_i(K\lalp)$ and 
the composite of $F\lsa\sub F_S$ with $P\lsa$ is
${\mathbf 1}_{F\lsa}$. By Lemma \ref{3.2},
$N_i\uan|_S=P_{S,\alpha}\sta(N\lsa\uan)$.
Now let $P_S$ be the projection
$$P_S= P\lalp\circ ({\mathbf 1}_{E\lst}\oplus P_{S,\alpha}):
G\lst\lra G\lsat\lra E\lst
$$
and let $Q_S=P_S\sta(M_i\ualg)$ be the pull back cone. 
Let $\tilde d_3$ be
\begin{equation*}
\begin{CD}
\tilde d_3: F_S|_w @>{P_{S,\alpha}|_w}>> F\lsa|_w @>{d_3}>> \cT_S^2|_w,
\end{CD}
\end{equation*}
then clearly we have a commutative diagram of vector spaces
\begin{equation}
\begin{CD}
F_S|_w @>>> G\lst|_w @<<< E\lst|_w\\
@VV{\tilde d_3}V @VV{P_S|_w}V @VV{d_1}V\\
\cT^2_S|_w @= \cT^2_S|_w @= \cT^2_S|_w.
\end{CD}
\end{equation}
Because $w$ is arbitrary, similar to the previous case, we have that $F_S$ 
intersects $Q_S$ transversally and $F_S\cap Q_S=N_i\uan|_S$, as
stratifiable currents.

To enable us to patch $Q_S$, where $S\sub Z_i$, to form a current in $G_{i,2}=
E_{i,2}\oplus F_i$, we need to show that $Q_S$ is independent of the
choice of analytic chart $\alpha$. Namely if we let $\beta\in\Xi$ be another
analytic chart so that $\iota_i(S)\sub\iota_{\beta}(Z_{\beta})$, then the cone current
$Q_S\pri\sub G\lst$ constructed using $F_{\beta}$, etc., is identical to $Q_S$.
Again, following the same argument before, it suffices to show that
the homomorphism $\tilde d_3\mh F_S|_w\to\cT^2_S|_w$ does
not depend on the choice of $\alpha$. Note that $\tilde d_3$ also fits into
the commutative diagram of exact sequences
\begin{equation}
\begin{CD}
T_{\rho\lalp(w)}R\lalp @>{d\phi\lalp(\rho\lalp(w))}>> F\lalp|_{\rho\lalp(w)}
@>>> \cT\lalp^2|_{\rho\lalp(w)}@>>> 0\\
@VVV @VVV @| \\
T_w R_i @>{d\phi_i(w)}>> F_S|_w @>>> \cT_S^2|_w @>>> 0.
\end{CD}
\lab{eq:3.33}
\end{equation}
Now assume  $\beta\in\Xi$ as before. Without loss of generality, 
we can assume that 
near $w$, the vector subbundles $\rho\lalp\sta F\lalp$ and 
$\rho_{\beta}\sta F_{\beta}$ span 
a $2l$-dimensional subvector bundle of $F_i$. 
Now let $V\lalp\to \tilde U\lalp$ and $V\lbe\to\tilde U\lbe$ 
be the vector bundles that define
$R\lalp$ and $R\lbe$ as in section 2 and let $V\lalpbe\to \tilde U_i$ be the
direct sum of the pull back of $V\lalp$ and $V\lbe$ via the tautological map 
$\tilde U_i\to \tilde U\lalp$ and $\tilde U_i\to\tilde U\lbe$. Then near
a neighborhood of $w\in\tilde U_i$, the set 
$\tilde\Phi\upmo(V\lalpbe)$ will form a base
of a smooth approximation containing $w$. We denote
$R\lalpbe=\tilde\Phi_i\upmo(V\lalpbe)$ and let 
$\phi\lalpbe\mh R\lalpbe\to V\lalpbe|_{R\lalpbe}$ be the lift of
$\tilde\Phi_i$.
Clearly, $R_i$ is still finer than $R\lalpbe$. Hence we have commutative
diagrams
\begin{equation}
\begin{CD}
T_{\rho\lalp(w)}R\lalp @>{d\phi\lalp(\rho\lalp(w))}>> V\lalp|_{w}
 @>>> \cT^2\lalp|_{\rho\lalp(w)} @>>> 0\\
@VVV @VVV @|\\
T_wR\lalpbe @>{d\phi\lalpbe(w)}>> V\lalpbe|_w @>>> \cT^2_i|_{w}
@>>> 0\\
@VVV@VVV@|\\
T_w R_i @>{d\phi_i(w)}>> F_i|_w @>>> \cT^2_i|_w @>>> 0\\
\end{CD}
\lab{eq:3.34}
\end{equation}
with exact rows.
Note that $V\lalpbe|_w\to \cT^2_i|_w$ is equal to
$$V\lalp|_{\rho\lalp(w)}\oplus V\lbe|_{\rho\lbe(w)}\lra
\Gamma(\om\lcpt^{0,1}(f_w\sta T_X))\lra H^{0,1}_{\dbar}(f_w\sta T_X)
\lra \cT^2_i|_w.
$$
(Here that $V\lalpbe|_w\to \cT_i^2|_w$ is defined apriori but 
not $F_i|_w\to\cT^2_i|_w$ because elements of $V\lalp|_w$
and $V\lbe|_w$ are $(0,1)$-forms with compact support.)
Therefore, the homomorphism $\tilde d_3$ defined earlier is 
independent of the choice of $\alpha$.

Now we are ready to prove the theorem. Let $i\in\Lambda$ be any 
approximation and let $\{S_a\}$ be an open covering 
of $Z_i$ so that to each $a$ there is an
$\alpha_a\in\Xi$ so that
$\iota_i(S_a)\sub\iota_{\alpha_a}(Z_{\alpha_a})$. We let
$G_{i,2}=E_{i,2}\oplus F_i$ and let 
$Q_{S_a}$ be the cone in
$G_{i,2}|_{S_a}$ constructed before
using the analytic chart $\alpha$. We know that over 
$G_{i,2}|_{S_a\cap
S_b}$, the currents $Q_{S_a}$ and $Q_{S_b}$ coincide. 
Hence $\{Q_{S_a}\}$ patchs
together to form a stratifiable current, denoted $Q_i$. Assume that
$j<i\in\Lambda$ be any two indices. Let $Z_{i,j}\sub Z_i$ be the open
subset $\iota_i\upmo(\iota_j(Z_j))$ and let $f^j_i\mh Z_{i,j}\to Z_j$
be the map induced by $Z_i$ being finer than $Z_j$. Then
$(f_i^j)\sta(F_j)$ is canonically a subbundle of $F_i|_{Z_{i,j}}$,
and $(f_i^j)\sta(E_{j,2})$ is 
canonically isomorphic to $E_{i,2}|_{Z_{i,j}}$. Let 
$(f_i^j)\sta(G_{j,2})\to
G_{i,2}|_{Z_{i,j}}$ be the induced homomorphism. It follows from the 
previous argument that $Q_i$
intersects $(f_i^j)\sta(G_{j,2})$ transversally and the intersection
$Q_i\cap(f_i^j)\sta(G_{j,2})$ is $(f_i^j)\sta(Q_j)$. Finally, by our
construction, $Q_i$ intersects transversally with $E_{i,2}$ and
$F_i\sub G_{i,2}$, and $E_{i,2}\cap Q_i=M_i\ualg$ and $F_i\cap
G_i=N_i\uan$. Let $\cG$ be the semi-$\QQ$-vector bundle $\{G_{i,2}\}$,
which is $\cE\oplus\cF$, and let $\cQ$ be the cone $\{Q_i\}$. 
It follows from the
perturbation argument in section two that for generic sections
$\eta_E$, $\eta_F$ and $\eta_G$ of $\cE$, $\cF$ and $\cG$ 
respectively, we have
$$[\mgna]\vir=[\eta_E\sta\cM\ualg]=[\eta_G\sta\cQ]
=[\eta_F\sta\cN\uan]=e[\Phi\mh\bB\to\bE].
$$
This proves the comparison theorem.

\end{document}